\documentclass{amsart}
\usepackage{amsmath,amsthm,amssymb,graphicx,mathtools,cancel,mathrsfs}
\usepackage[pdfusetitle,breaklinks]{hyperref}
\usepackage{marginnote,bm}
\usepackage[shortlabels]{enumitem}
\usepackage{xcolor}
\usepackage[capitalize]{cleveref}
\usepackage{microtype}

\newcommand{\R}{\mathbb{R}}
\renewcommand{\P}{\mathbb{P}}
\newcommand{\Z}{\mathbb{Z}}
\newcommand{\Q}{\mathbb{Q}}

\newcommand{\F}{\mathbb{F}}
\newcommand{\A}{\mathbb{A}}

\newcommand{\pp}{\mathfrak{p}}

\newcommand{\km}{\mathfrak{m}}

\newcommand{\kp}{\pp}

\newcommand{\OO}{\mathcal{O}}

\newcommand{\cI}{\mathcal{I}}
\newcommand{\cT}{\mathcal{T}}
\newcommand{\cB}{\mathcal{B}}
\newcommand{\cG}{\mathcal{G}}

\newcommand{\bt}{\mathbf{t}}
\newcommand{\bG}{\mathbf{G}}

\newcommand{\bC}{\mathbf{C}}

\newcommand{\bSig}{\mathbf{\Sigma}}

\newcommand{\wh}{\widehat}

\DeclareMathOperator{\lcm}{lcm}

\DeclareMathOperator{\Tr}{Tr}
\DeclareMathOperator{\id}{id}
\DeclareMathOperator{\ord}{ord}

\DeclareMathOperator{\Mat}{Mat}

\DeclareMathOperator{\sgn}{sgn}

\DeclareMathOperator{\NP}{NP}
\DeclareMathOperator{\GNP}{GNP}
\DeclareMathOperator{\GL}{GL}

\DeclareMathOperator{\Sol}{Sol}
\DeclareMathOperator{\Sym}{Sym}

\DeclareRobustCommand{\littlep}{\ensuremath{p}}
\AddToNoCaseChangeList{\littlep}

\theoremstyle{plain}
\newtheorem{theorem}{Theorem}[section]
\newtheorem{corollary}[theorem]{Corollary}
\newtheorem{lemma}[theorem]{Lemma}

\newtheorem{proposition}[theorem]{Proposition}

\theoremstyle{definition}
\newtheorem{definition}[theorem]{Definition}
\newtheorem{example}[theorem]{Example}
\newtheorem{remark}[theorem]{Remark}
\newtheorem*{remark*}{Remark}
\newtheorem*{example*}{Example}

\AddToHook{env/corollary/begin}{\crefalias{theorem}{corollary}}
\AddToHook{env/lemma/begin}{\crefalias{theorem}{lemma}}
\AddToHook{env/keylemma/begin}{\crefalias{theorem}{keylemma}}
\AddToHook{env/proposition/begin}{\crefalias{theorem}{proposition}}
\AddToHook{env/claim/begin}{\crefalias{theorem}{claim}}
\AddToHook{env/definition/begin}{\crefalias{theorem}{definition}}
\AddToHook{env/example/begin}{\crefalias{theorem}{example}}
\AddToHook{env/remark/begin}{\crefalias{theorem}{remark}}

\title[Generic elementary abelian 
\texorpdfstring{\littlep}{p}-covers]
{A global generic locus for elementary abelian 
    \texorpdfstring{\littlep}{p}-covers of 
    \texorpdfstring{$\P^1$}{P1} with 
two ramification points}

\author{Robert Moore}
\author{Hui June Zhu}
\date{\today}

\address{
Robert Moore,
Department of mathematics,
State University of New York at Buffalo,
Buffalo, NY 14260. 
The United States.
}
\email{rcmoore@buffalo.edu}

\address{
Hui June Zhu,
Department of mathematics,
State University of New York at Buffalo,
Buffalo, NY 14260. 
The United States.
}
\email{hjzhu@math.buffalo.edu}

\begin{document}
\begin{abstract}
    In this paper, we study elementary abelian $p$-covers
    of $\P^1$ totally ramified at exactly two points and unramified elsewhere.
    Such covers are isomorphic to the rank-$\ell$ Artin--Schreier curves
    \[
        X_{f, \ell} : y^{p^\ell} - y = f(x), 
    \]
    where $f \in \overline{\F}_p[x^{\pm1}]$
 has a pole of order $d_1>0$ at $0$ 
and order $d_2>0$ at $\infty$.
    After fixing $d_1$ and $d_2$,
    the cover $X_{f, \ell}$ is generic
    if its Newton polygon is minimal
    as $f$ varies over all Laurent
    polynomials with pole orders $d_1$ and $d_2$.
    For $p$ sufficiently large relative to $d_1$ and $d_2$,
    we construct a combinatorial cost polygon $\GL(d_1,d_2,p)$ 
    and prove that its attainment is characterized by the nonvanishing 
    of a homogeneous polynomial in the coefficients of $f$.
    Moreover, we prove that these conditions stabilize:
    there exists a Zariski-open dense subset $\mathcal{U}$
    in the coefficient space that is defined over $\Q$
    and depends only on $d_1$ and $d_2$,
    whereby for any $f$ whose coefficients lie in
    $\mathcal{U}(\overline{\Q})$,
    the reduction is generic at every prime ideal above sufficiently large 
    rational primes, and for every rank $\ell$ permitted by the residue field.
    This result extends Zhu's result for one ramification point 
    in \cite{zhu2025}.
\end{abstract}

\maketitle

\section{Introduction}
Let $q = p^b$ for a rational prime number $p$ and a positive integer $b$.
If $C$ is a smooth projective curve over the finite field $\F_q$ 
with $q$ elements, then the zeta function $Z(C / \F_q, s)$
is classically defined as the generating function
$\exp \left(\sum_{m=1}^\infty \#C(\F_{q^m}) \frac{s^m}{m}\right)$.
By the theorems of Weil (see \cite{ref:weil}),
the zeta function is a rational function of the form
\[
    \frac{P_C(s)}{(1-s)(1-qs)}, \quad P_C(s) = 
    \prod_{i=1}^{2g} (1-\alpha_i s) \in 1 + s\Z[s],
\]
where $g$ is the genus of $C$, and the reciprocal roots $\alpha_i$
have complex absolute value $\sqrt{q}$.

The $q$-adic absolute value of the reciprocal roots $\alpha_i$
admits a well-known symmetry.
Let $P_C(s) = 1 + c_1 s + \dots + c_{2g} s^{2g}$.
The ($q$-adic) \textbf{Newton polygon} of this polynomial is
the lower convex hull
of the set $\{(i, v_q(c_i)) \mid i = 0, \dots, 2g\}$,
where we set $c_0 = 1$,
and $v_q$ denotes the $q$-adic valuation on $\overline{\Q}_p$,
normalized so that $v_q(q) = 1$.
It is a piecewise linear function, with endpoints at $(0,0)$ and $(2g, g)$,
and it is symmetric in the sense that for every segment
with horizontal length $h$ and slope $\lambda$,
there is a corresponding segment with horizontal length $h$
and slope $1 - \lambda$.
We denote by $\NP(C / \F_q)$ the 
$q$-adic Newton polygon of the numerator of the zeta
function of $C$, and simply refer to it as the \textit{Newton polygon of $C$}.
The Newton polygon for $C$ is moreover invariant under choice of
the base field $\F_q$, provided the $q$-adic valuation is adjusted accordingly.
Thus we may simply write $\NP(C)$ for the Newton polygon of $C$,
with the base field being any finite field compatible with $C$.

Let $f \in \overline{\F}_p(x)$ 
be a rational function,
and let $\ell \ge 1$.
We define the \textit{rank-$\ell$ Artin--Schreier curve} $X_{f, \ell}$ to be
the smooth projective curve defined by the affine equation 
\[
    X_{f, \ell} : y^{p^\ell} - y = f.
\]
All elementary abelian $p$-covers of $\P^1$ are represented by equations of this
form (see \cite{ref:garcia1991}),
making them a subject of geometric interest.

This paper is concerned with elementary abelian $p$-covers
of $\P^1$ with two ramification points and unramified elsewhere ---
which, after applying a suitable isomorphism
sending the ramification points to $0$ and $\infty$,
can be characterized by
 $X_{f, \ell}$ with $f \in \overline{\F}_p[x^{\pm1}]$.
Let $f(x) = \sum_{i=-d_1}^{d_2} a_i x^i$, for $d_1, d_2 \ge 1$
and $a_{-d_1}a_{d_2}\ne 0$.
We also assume that $p \nmid d_1 d_2$.
Then there is a pole of order $d_1$ at $0$ and a pole of order $d_2$
at $\infty$.
Its genus equal $g=\frac{(p^\ell-1)(d_1+d_2)}{2}$.
Let $\mathscr{X}_{\overline{\F}_p, \ell}(d_1, d_2)$ denote the set of all
such rank-$\ell$ Artin--Schreier curves over $\overline{\F}_p$.
Then we define the \textbf{generic Newton polygon}
\begin{equation}
    \label{E:GNP}
    \mathrm{GNP}_{\overline{\F}_p, \ell}(d_1, d_2)  \coloneqq 
    \inf_{X \in \mathscr{X}_{\overline{\F}_p, \ell}(d_1, d_2)} \NP(X).
\end{equation}
Here, the infimum is taken with respect to the standard 
pointwise partial ordering on the set $(\R_{\ge 0})^{[0,2g]}$
of functions $h : [0,2g] \to \R_{\ge 0}$.
Grothendieck's specialization theorem implies that
there exists a curve $X$ that achieves this infimum.
If $\NP(X) = \mathrm{GNP}_{\overline{\F}_p, \ell}(d_1, d_2)$,
then we say that $X$ is \textbf{generic}.

The subject of determining the shape of 
$\GNP_{\overline{\F}_p, \ell}(d_1, d_2)$,
as well as the set of generic curves in 
$\mathscr{X}_{\overline{\F}_p, \ell}(d_1, d_2)$,
has been studied extensively for the past several decades.
The cohomological methods of Adolphson and Sperber in \cite{AS89}
yield a universal lower bound known as the Hodge polygon.
However, this lower bound is only sharp for very specific cases.

As $d_1$, $d_2$, and $p$ vary,
$\mathrm{GNP}_{\overline{\F}_p, \ell}(d_1, d_2)$
varies in ways that have proven difficult to generalize.
However, when $d_1$ and $d_2$ are fixed,
the asymptotic behavior
of $\mathrm{GNP}_{\overline{\F}_p, \ell}(d_1, d_2)$
stabilizes
as $p \to \infty$.
For the rank $1$ case (i.e. $\ell = 1$)
for polynomials,
Zhu
constructed the generic Newton polygon in 
\cite{zhu2002}
for all $p$ sufficiently large,
and defined a global open condition
that characterizes the locus of generic Artin--Schreier curves.
Blache and F\'erard refined this result 
in \cite{ref:blache2007}.

Most recently, Zhu discovered the locus of generic elementary abelian
$p$-covers with one ramification point in \cite{zhu2025}.
Consider a polynomial $f \in \overline{\Q}[x]$.
Let $\kp$ be a prime ideal in some number field $K$ containing the coefficients
of $f$.
If the coefficients of $f$ are $\kp$-integral,
then one can consider the image $\overline{f}$
of $f$ in the residue field $k_{\kp} = \OO_K / \kp$,
and determine whether the curve $X_{\overline{f}, \ell}$
is generic.
Zhu proved in \cite[Theorem 1.2]{zhu2025} that an open condition
on polynomials in $\overline{\Q}[x]$ defines a locus of polynomials
whose reductions yield generic curves for \textit{every}
$\kp$ above a sufficiently large $p$.

Our first main theorem extends Zhu's global result on elementary abelian $p$-covers 
with one ramification point to those with two
ramification points.
\begin{theorem}
    \label{T:ordinary-curves}
    Let $d_1, d_2 \ge 1$.  Let $I = \{-d_1, \dots, d_2\}$,
    and let $\A^I_{\Q}$ denote affine $\Q$-space indexed by
    $I$.
    For any $a \in \A^I_{\Q}(\overline{\Q})$,
    let $f_a = \sum_{i \in I} a_i x^i$.
    Then there exists a homogeneous Zariski-open and dense subset 
    $\mathcal{U} \subseteq \A_{\Q}^I$ defined over $\Q$
    for which the following are equivalent
    for every $a \in \A_\Q^I(\overline{\Q})$
    with $a_{-d_1}a_{d_2} \ne 0$:
    \begin{enumerate}[(1)]
        \item $a \in \mathcal{U}(\overline{\Q})$,
        \item There is a constant $N_a$ such that,
            for every prime number $p > N_a$,
            every number field $K / \Q(a_i \mid i \in I)$,
            and every prime ideal $\kp$ in $K$
            above $p$, the reduction $\overline{f}_a = f_a \pmod{\kp}$
            is defined, and the curve $X_{\overline{f}_a, \ell}$ is generic for
            every $\ell$ that divides $[k_\kp : \F_p]$.
    \end{enumerate}
\end{theorem}
This extension cannot be
obtained merely by treating positive and negative
terms of $f$ independently and reducing to the one-pole case on each sign.
In fact, 
the open set $\mathcal{U}$
generally does not factor into open sets independently 
defined on the positive and negative indices, because
the positive and negative terms interact with one another
in relevant ways that generally cannot be separated.
This is illustrated in \cref{R:cross-sign-minimizers}.

We construct a combinatorial cost polygon $\GL(d_1,d_2,p)$
in Definition \ref{D:ANP}, so that for $p$ large enough
$$
\GNP_{\bar{F}_p,\ell}(d_1,d_2)(x)=(p^\ell-1)\GL(d_1,d_2,p)
\left(\frac{x}{p^\ell-1}\right).
$$

The same global open condition also controls $T$-adic exponential sums.
Combining our construction of $\mathcal U$ with the interpolation argument
of \cite{ref:liu2016} extended in \Cref{P:lln-result} 
to the range $p\ge2D$, gives the following consequence, 
proved as \Cref{T:T-adic}.
We use the $C$- and $L$-functions defined in
\cref{SS:T-adic-intro} 
and write $\pi_m=\zeta_{p^m}-1$, where $\zeta_{p^m}$ is a
primitive $p^m$-th root of unity.

\begin{corollary}[Simultaneous $T$-adic genericity]
\label{cor:introduction-T-adic}
Let $I$, $f_a$, and $\mathcal{U}$ be as in \cref{T:ordinary-curves}, and let
$a\in\A^I_{\Q}(\overline{\Q})$ with
$a_{-d_1}a_{d_2}\ne 0$.
Then $a\in \mathcal{U}(\overline{\Q})$ if and only
if there is a constant $N_a$ such that, for every prime $p>N_a$, every
number field $K$ containing the coefficients of $f_a$, and every prime
ideal $\kp$ of $K$ above $p$, the reduction $\bar f_a$ is defined,
and the functions
\[
  C^*(\bar f_a;T,s),\qquad
  \left.C^*(\bar f_a;T,s)\right|_{T=\pi_m},\qquad
  \left.L^*(\bar f_a;T,s)\right|_{T=\pi_m}
\]
attain their respective generic Newton polygons for every $m\ge1$.
Here the first polygon is taken with respect to the $T$-adic valuation,
and the specialized polygons with respect to the $\pi_m$-adic valuation,
with the pole orders $d_1$ and $d_2$ fixed. The polygon of the specialized
$L$-function is considered on $[0,p^{m-1}(d_1+d_2)]$.
\end{corollary}

In particular, the bound $N_a$ is independent of the conductor level $m$
and of the number field $K$. Geometrically, the specialized $L$-function
is the factor associated with a character of order $p^m$ in the cyclic
Artin--Schreier--Witt tower determined by the coefficientwise
Teichm\"uller lift of $\bar f_a$. This connects our result with the study
of slopes in Artin--Schreier--Witt towers. Ren, Wan, Xiao, and
Yu~\cite{ref:RWXY18} study Newton slopes and spectral halos for
higher-rank Artin--Schreier--Witt towers. Tang~\cite{ref:tang2022} proves
slope uniformity for strongly genus-stable higher-rank towers and slope
stability under an additional coefficient-growth condition. These works
study character variation in fixed characteristic for towers ramified
only at infinity. Our
$T$-adic result complements this theory by giving a single open locus
over $\Q$ for two-pole Laurent polynomials whose reductions are
generic at every finite conductor level for all sufficiently large
primes.

Our proof is based on valuation estimates for Dwork Fredholm
determinants. \cref{S:appendix} establishes the combinatorial results that
identify the minimal terms and prevent their cancellation.
\cref{S:dwork-work}
reviews the relevant Dwork theory and $T$-adic interpolation.
\cref{S:ck_formula}
constructs the polynomials controlling genericity and extracts
factors independent of $p$, while
\cref{S:extension_theorem} proves the compatibility of
the genericity criterion with finite extensions of the base field.
\cref{S:independence} combines these results to prove
\cref{T:ordinary-curves}. After establishing in \cref{P:lln-result} simultaneous genericity for
the $T$-adic $C$-function and all the specialized $C$- and $L$-functions using the method of \cite{ref:liu2016},
we prove \Cref{cor:introduction-T-adic} in \Cref{T:T-adic}.

\section{Combinatorics of minimal terms in Fredholm determinants}
\label{S:appendix}
In this section, we prove the combinatorial facts needed for
analyzing the Fredholm determinants in the sections that follow.
We first describe integer representations and minimal representations
in \cref{SS:cmp}.
All of the key results are in \cref{SS:cycle-sums} and \cref{SS:ordering}.
\cref{SS:cycle-sums} will determine which index sets and
permutations contribute to the valuations in the Fredholm determinant.
It will also determine which points
are vertices on the generic Newton polygon.
Finally, we will choose a monomial order
in \cref{SS:ordering}
that singles out a noncancelling term in the relevant part of the
Fredholm determinant's coefficients.

For any integer $N$
and any positive integer $d$,
we let $[N]_d$ denote the remainder of $N$
modulo $d$ --- 
that is, $[N]_d$ is the unique element 
in $\{0, \dots, d-1\}$ that is congruent to $N$
modulo $d$.
For the duration of this paper, we will fix
positive integers $d_1$ and $d_2$,
and set
\[
    D = \lcm(d_1,d_2),\quad 
    I = \{-d_1, \dots, d_2\}, \quad \mbox{and} \quad 
    I_0 = I \setminus \{0\}.
\]

\subsection{Minimal integer representations}
\label{SS:cmp}
The Dwork operator used in \cref{S:dwork-work,S:ck_formula}
is built out of a $p$-adic Laurent
power series $F(s) = \sum_{N \in \Z} G_N s^N$,
and the coefficients $G_N$
are obtained by collecting the ways to represent $N$
as a nonnegative integer linear combination of
the support elements.
We will be interested in minimizing that linear combination,
which amounts to a change-making problem.
We will describe that problem, and its solution,
in this subsection.

In computational mathematics,
a change-making problem concerns the minimization
of a summation $\sum_i t_i$
subject to the constraint that 
$t_i$ are nonnegative integers,
and $\sum_i t_i c_i = N$
for some coin set $\{c_i\}$ of positive integers,
and some fixed nonnegative integer $N$.
In this section, we will consider an ``extended''
change-making problem, where the coins $c_i$
are both positive and negative, and $N$ is any integer.

\begin{definition}
    \label{D:cmp}
    For $N \in \Z$,
    define
    \[
        \cT(N) = \left\{\bt = (t_i) \in \Z_{\ge 0}^{I_0} \mid 
        \sum_{i \in I_0} i t_i = N\right\}, \quad 
        |\bt| = \sum_{i \in I_0} t_i, \quad 
        M(N) = \min_{\bt \in \cT(N)} |\bt|,
    \]
    and $\Sol(N) = \left\{\bt \in \cT(N) \mid |\bt| = M(N)\right\}$.
    An element $\bt \in \cT(N)$
    is called a \textbf{representation of $N$}.
    If $\bt \in \Sol(N)$,
    then $\bt$ is a \textbf{solution} for $N$,
    or a \textbf{minimal representation} of $N$.
\end{definition}

We first prove two lemmas about $M(N)$
and $\Sol(N)$.
\begin{lemma}
    \label{L:cmp-facts}
    For any $N \in \Z$, set $d = d_1$ if $N<0$ 
    and $d = d_2$ if $N \ge 0$. 
    Then $M(N) = \left\lceil |N|/d \right\rceil$.
    Moreover, for nonzero $N$,
    every $\bt \in \Sol(N)$ has $t_i = 0$
    for all $i$ for which $\sgn(i) \ne \sgn(N)$.

    There is a distinguished solution $\bt^+ = (t^+_i)_{i \in I_0}$
    of $N$,
    given as follows:
    \begin{itemize}
        \item If $N \ge 0$, then $t^+_{d_2} = \lfloor N/d_2\rfloor$,
            and $t^+_{[N]_{d_2}} = 1$ if $[N]_{d_2} \ne 0$.
            $t^+_j = 0$ for all other $j$.
        \item If $N < 0$, then $t^+_{-d_1} = \lfloor -N / d_1\rfloor$,
            and $t^+_{-[-N]_{d_1}}=1$ if $[-N]_{d_1} \ne 0$.
            $t^+_j = 0$ for all other $j$.
    \end{itemize}
\end{lemma}
\begin{proof}
    Suppose $N\ge0$, and let $\bt \in \Sol(N)$.
    Write $T_+ = \sum_{i > 0} t_i$
    and $T_- = \sum_{i < 0} t_i$.
    Since $N \le d_2 T_+$, we have that
    \[
        |\bt| = T_+ + T_- \ge T_+ \ge \left\lceil N/d_2\right\rceil.
    \]
    The division algorithm gives the representation $\bt^+$ for $N$
    with this cardinality, so indeed
    $|\bt| = \left\lceil N / d_2 \right\rceil$.
    Moreover, we must have $T_- = 0$ since otherwise
    $|\bt| = T_+ + T_- > T_+ \ge \left\lceil N / d_2\right\rceil$.
    The proof for $N<0$ is symmetric.
\end{proof}

\begin{lemma}
    \label{L:lower-term-bijection}
    For nonnegative integers $d$ and $n$,
    let $\cT_d(n)$ denote the set of nonnegative integer tuples
    $(t_1, \dots, t_d) \in \Z_{\ge 0}^d$ for which 
    $\sum_{i=1}^d t_i i = n$,
    with $\cT_0(0) = \{()\}$ and $\cT_0(n) = \varnothing$ for positive $n$.

    Let $N \in \Z \setminus\{0\}$ and set 
    $d = d_1$ or $d = d_2$, depending on whether $N$ is negative or positive,
    respectively.
    Let $\Delta = d M(N) - |N|$,
    and suppose $M(N) \ge d-1$.
    Then $\Sol(N)$ and $\cT_{d-1}(\Delta)$
    are in bijection.
    $\Sol(N) \to \cT_{d-1}(\Delta)$ is given by 
    $(t_i)_{i \in I_0} \mapsto (t_{\sgn(N)(d-i)})_{i=1}^{d-1}$,
    the inverse map from $\mathbf{u} = (u_j) \in \cT_{d-1}(\Delta)$
    is given by
    $t_{\sgn(N)j} = u_{d-j}$ for $j=1,\dots,d-1$,
    $t_{\sgn(N)d} = M(N) - |\mathbf{u}|$.
\end{lemma}
\begin{proof}
    Every minimal representation of $N$ has size $M(N)$
    and uses indices of sign $\sgn(N)$ by \cref{L:cmp-facts}.
    Therefore 
    \[
        \Delta = d M(N) - |N| = \sum_{i=1}^{d-1} (t_{\sgn(N)(d-i)})i.
    \]
    The stated formulas give the inverse map.
\end{proof}

\subsection{Permutation costs and minimizing index sets}
\label{SS:cycle-sums}
Evaluating the valuation of the Fredholm determinant's coefficients 
concerns the determinants of finite principal minors of
the infinite matrix $(G_{pi-j})_{i,j \in \Z}$,
where the $p$-adic valuation of $G_N$ is generally 
proportional to the solution value $M(N)$ described in the previous
subsection.
Approximating the $p$-adic valuation of these determinants, therefore,
amounts to minimizing summations of the form
$\sum_{\ell \in \cI} M(p\ell-\sigma(\ell))$,
where $\cI$ is a subset of $\Z$ with fixed finite size
and $\sigma$ is a permutation on $\cI$.
This subsection analyzes that problem.
\cref{P:small-one} provides a concrete solution when $\cI$
is an interval, and \cref{P:cI-is-linspace}
demonstrates that for sufficiently large $p$,
an interval will be the minimizing index set $\cI$.
Then \cref{C:kplus1} establishes the convexity of the generic construction,
while \cref{P:minimal-bijections} will ultimately validate a reduction 
we employ in \cref{S:extension_theorem}.

We inherit the notation of \cref{D:cmp}
from the previous subsection.

\begin{definition}
    \label{D:cycle-sums}
    Let $p$ be a prime number.
    For a finite subset $\cI$ of $\Z$,
    and a permutation $\sigma$ on $\cI$,
    we define 
    \[
        S^p(\cI) \coloneqq \min_{\sigma \in \mathrm{Sym}(\cI)} 
        \sum_{\ell \in \cI} M(p\ell - \sigma(\ell)).
    \]
    For a nonnegative integer $k$,
    we set 
    \[
        S^p(k) \coloneqq \min_{\cI \in \binom{\Z}{k}} S^p(\cI),
    \]
    and 
    \[
        \bSig^p(\cI) = \left\{\sigma \in \mathrm{Sym}(\cI) \mid
        \sum_{\ell \in \cI}M(p\ell - \sigma(\ell)) = S^p(\cI)\right\},
    \]
    where $\binom{\Z}{k}$ denotes the set of all subsets of $\Z$
    with cardinality $k$.
    We will often write $\sigma\ell$ instead of $\sigma(\ell)$
    for ease of reading.
\end{definition}

The function $S^p(k)$ defines, up to scale,
the $y$-coordinates of the combinatorial cost polygon,
which we also define now.
\begin{definition}[Combinatorial cost polygon]
    \label{D:ANP}
    For a prime number $p$ and integers $d_1, d_2 \ge 1$,
    we define $\GL(d_1,d_2,p)$
    to be the lower convex hull of the 
    set
    \[
        \left\{ \left(k, \frac{S^p(k)}{p-1}\right) \mid 
        k \in \{1, \dots, d_1+d_2\}\right\} \cup \{(0,0)\}.
    \]
    We call this the \textit{combinatorial cost polygon}
    of $(d_1, d_2)$ on $p$.
    
    When $d_1$ and $d_2$ are understood, we 
    define $V_p \subseteq \{1, \dots, d_1+d_2 - 1\}$
    to be the set of abscissae of internal vertices of the combinatorial cost polygon 
    $\GL(d_1,d_2,p)$.
\end{definition}
One of the goals of this subsection is to 
prove the convexity of the function $S^p(k)$,
and to describe an equivalent condition for an element $k$ to be in $V_p$.
Both tasks are accomplished in \cref{C:kplus1}.

We begin our exploration by deriving a more computable 
expression for $S^p(\cI)$ when $\cI$ is a sufficiently narrow interval.
\begin{proposition}
    \label{P:small-one}
    Suppose $p$ is coprime to both $d_1$ and $d_2$,
    and let $k_1$ and $k_2$ be nonnegative integers
    such that 
    $-d_1 \le -k_1 \le 0 \le k_2 \le d_2$.
    If $\cI = [-k_1, k_2] \cap \Z$, then
    \[
        S^p(\cI) = \sum_{\ell=1}^{k_1} 
        \left\lceil \frac{p\ell - k_1}{d_1} \right\rceil + 
        \sum_{\ell = 1}^{k_2} \left\lceil \frac{p\ell - k_2}{d_2} \right\rceil.
    \]
    Moreover, if $\sigma \in \bSig^p(\cI)$,
    then $\sigma(0) = 0$, and
    $\left\lceil \frac{p\ell - \sigma\ell}{d_2}\right\rceil = 
    \left\lceil \frac{p\ell - k_2}{d_2}\right\rceil$
    for all positive $\ell$,
    and $\left\lceil \frac{p(-\ell) + \sigma\ell}{d_1}\right\rceil
    = \left\lceil \frac{p(-\ell) - k_1}{d_1}\right\rceil$
    for all negative $\ell$.
\end{proposition}
\begin{proof}
    Note that the righthand side is certainly a lower bound
    for $S^p(\cI)$,
    since for every positive $\ell$,
    $p\ell - \sigma\ell \ge p\ell - k_2$,
    and so $M(p\ell - \sigma\ell) \ge \left\lceil
        \frac{p\ell-k_2}{d_2}\right\rceil$.
    Similarly, for every negative $\ell$,
    $M(p\ell - \sigma\ell) \ge \left\lceil
    \frac{p(-\ell)-k_1}{d_1}\right\rceil$.
    And for $\ell = 0$, $M(p\ell - \sigma\ell) \ge 0$.
    It remains to show that this lower bound is achieved for some permutation
    $\sigma$.

    We construct the permutation directly and refer to it as $\tau$.
    Starting with $\ell = k_2$ and working downward inductively to $\ell=1$,
    we let $I_\ell = \{1, \dots, k_2\} \setminus \{\tau^{-1}(j) \mid j \ge
    \ell+1\}$.
    Then we set $\tau^{-1}(\ell)$ equal to the $i \in I_\ell$
    that minimizes the value of $[pi-\ell]_{d_2}$.
    Since $p$ and $d_2$ are coprime,
    the set $\{pi-\ell \mid 1 \le i \le d_2\}$ forms a complete residue system
    modulo $d_2$; and since $|I_\ell| \le k_2 \le d_2$,
    $[pi-\ell]_{d_2}$ will be minimized for a unique choice of $i$.
    Thus $\tau^{-1}$ is well-defined and injective.

    Suppose that for some $i$, 
    $\left\lceil \frac{pi - \tau(i)}{d_2}\right\rceil >
    \left\lceil \frac{pi - k_2}{d_2}\right\rceil$.
    Then for some $\tau (i) + 1 \le \ell \le k_2$,
    $pi-\ell \equiv 0 \pmod{d_2}$,
    and so $[pi-\ell]_{d_2} = 0$.
    Since $\ell \ne \tau(i)$,
    we must have that $i \not\in I_\ell$,
    hence $i \in \{\tau^{-1}(j) \mid j \ge \ell+1\}$,
    hence $i = \tau^{-1}(j)$ for some $j \ge \ell+1$,
    hence $\tau(i) \ge \ell+1$.
    But then $\ell+1 \le \tau(i) < \tau(i)+1 \le \ell$,
    a contradiction.

    After performing a similar inductive construction for $\tau$ for
    negative $\ell$, and setting $\tau(0) = 0$,
    we obtain the desired result.
\end{proof}
\begin{remark}
    \label{R:SpI-mod-invariant}
    The necessary and sufficient conditions on 
    $\sigma$ for it to be in $\bSig^p(\cI)$
    in \cref{P:small-one}
    are modular conditions to $d_1$ and $d_2$.
    Consequently, if $p' \equiv p$ modulo $D$,
    then the same interval $\cI$ will satisfy 
    $\bSig^p(\cI) = \bSig^{p'}(\cI)$,
    assuming $p,p'>\max(k_1,k_2)$.
\end{remark}

Next we introduce a weight function, which will
aid in determining the minimizing index sets $\cI$.
\begin{definition}
    \label{D:weight-function}
    We define a weight function $w : \Z \to \Q_{\ge 0}$ by
    \[
        w(i) = 
        \begin{cases*}
            \frac{-i}{d_1} & if $i \le 0$ \\
            \frac{i}{d_2} & if $i > 0$
        \end{cases*}
    \]
    and extend it to finite subsets $\cI \subseteq \Z$
    by $w(\cI) \coloneqq \sum_{\ell \in \cI} w(\ell)$.
    For every $k \ge 1$,
    we define
    $W_k \coloneqq \min_{\cI \in \binom{\Z}{k}} w(\cI)$.
    If a subset $\cI \in \binom{\Z}{k}$ 
    satisfies
    $w(\cI) = W_k$, then we say that $\cI$
    is \textit{weight-minimizing}.
\end{definition}

We provide some basic properties of our weight function.
\begin{lemma}
    \label{L:weight-properties}
    Let $k \ge 1$.
    \begin{enumerate}[(1)]
        \item For any $N \in \Z$, $M(N) = \lceil w(N)\rceil$.
        \item If $\cI \in \binom{\Z}{k}$ is weight-minimizing,
            then $\cI$ is an interval of the form
            \[
                \cI = \{-k_1, \dots, k_2\}
            \]
            for nonnegative integers $k_1$ and $k_2$
            that sum to $k-1$.
            The multiset $\{w(\ell) \mid \ell \in \cI \setminus\{0\}\}$
            consists of
            the $k-1$ smallest elements of 
            the multiset 
            $\{1/d_1, 2/d_1, \dots\} \sqcup \{1/d_2, 2/d_2, \dots\}$.
        \item There are either one or two distinct weight-minimizing
            sets in $\binom{\Z}{k}$.
        \item There are two distinct weight-minimizing sets 
            $\cI, \mathcal{J} \in \binom{\Z}{k}$ 
            if and only if there exist nonnegative integers $k_1$ and $k_2$
            for which $k_1 / d_1 = (k_2+1) / d_2$ and $k_1 + k_2 = k-1$.
            In this case, one of the subsets will be of the form 
            $\cI = \{-k_1, \dots, k_2\}$
            and the other will be of the form 
            $\mathcal{J} = \{-k_1+1, \dots, k_2+1\}$.
        \item $(p-1)w(\cI) \le S^p(\cI)$.
            In particular, $(p-1)W_k \le S^p(k)$.
        \item $S^p(k) < (p-1)W_k + k$.
    \end{enumerate}
\end{lemma}
\begin{proof}
    (1) - (3) are obvious.
    
    (4)
    Both of the ``if-and-only-if''
    conditions are equivalent to $\{1/d_1, 2/d_1, \dots\}$
    and $\{1/d_2, 2/d_2, \dots\}$ sharing a common element
    $r = k_1 / d_1 = (k_2+1) / d_2$ for which $k_1 + k_2 = k-1$.

    (5)
    Choose a minimizing permutation $\sigma$.
    Note that 
    \[
        M(p\ell - \sigma\ell) = \left\lceil w(p\ell - \sigma\ell)\right\rceil
        \ge w(p\ell - \sigma\ell) \ge p w(\ell) - w(\sigma\ell),
    \]
    so 
    \[
        S^p(\cI) = \sum_{\ell \in \cI} M(p\ell - \sigma\ell)
        \ge (p-1) \sum_{\ell \in \cI} w(\ell)
        = (p-1) w(\cI).
    \]
    If the choice of $\cI$ yields $S^p(\cI) = S^p(k)$,
    then we obtain $S^p(k) = S^p(\cI) \ge (p-1)w(\cI) \ge (p-1)W_k$.

    (6) Let $\cI' = \{-k_1, \dots, k_2\}$  be a weight-minimizing set,
    and take $\sigma'$ to be the identity.
    Then since
    $M(p\ell - \ell) = \lceil w(p\ell - \ell)\rceil
    < (p-1)w(\ell) + 1$, we have
    \[
        S^p(k) \le \sum_{\ell \in \cI'} M(p\ell -
        \ell) < (p-1) W_k + k.
    \]
\end{proof}

Next we introduce a key permutation lemma.
\begin{lemma}[Permutation lemma]
    \label{L:key}
    Let $1 \le k \le d_1 + d_2 + 1$ and $p > 2D$.
    Let $b \ge 1$ be an integer,
    and for every $m \in \Z / b\Z$,
    let $\cI_m \in \binom{\Z}{k}$.
    Construct the disjoint union of the $\cI_m$ as follows:
    \[
        \mathcal{A} = \bigsqcup_{m \in \Z / b\Z} \cI_m \times \{m\}.
    \]
    For any $a = (\ell, m) \in \mathcal{A}$, let $\pi_1(a) = \ell$.

    Then for any permutation $\tau : \mathcal{A} \to \mathcal{A}$,
    we have
    \[
        \sum_{a \in \mathcal{A}} M(p \pi_1(a) - \pi_1(\tau a)) \ge b S^p(k),
    \]
    with equality holding only if $\cI_m$
    is weight-minimizing for all $m \in \Z / b\Z$.
    In this case, we also have 
    \[
        \sum_{a \in \mathcal{A}} M(p \pi_1(a) - \pi_1(\tau a)) \ge 
        \sum_{m \in \Z / b\Z} S^p(\cI_m).
    \]
\end{lemma}
%% new proof

\begin{proof}
Write
$$
    C(\mathcal A,\tau)
    =\sum_{a\in\mathcal A}
      M\bigl(p\pi_1(a)-\pi_1(\tau a)\bigr),
    \qquad
    F_d(t)=\sum_{j=1}^{t}
      \left\lceil\frac{pj-t}{d}\right\rceil,
$$
with $F_d(0)=0$.

\smallskip
\noindent\emph{Weight-minimizing fibers.}
First suppose that every $\cI_m$ is weight-minimizing.
Since $k\le d_1+d_2+1$, each such set is contained in
$\{-d_1,\ldots,d_2\}$.

If all the fibers are equal to
$\cI=\{-k_1,\ldots,k_2\}$, then every target has first
coordinate in $\cI$. The row bounds in the proof of
\cref{P:small-one} therefore give
$$
    C(\mathcal A,\tau)
    \ge b\bigl(F_{d_1}(k_1)+F_{d_2}(k_2)\bigr)
    =bS^p(\cI)
    =\sum_m S^p(\cI_m)
    \ge bS^p(k).
$$
This includes the case $k=d_1+d_2+1$.

Otherwise, \cref{L:weight-properties} shows that the fibers
take the two values
$$
    \cI=\{-k_1,\ldots,k_2\},
    \qquad
    \mathcal J=\{-k_1+1,\ldots,k_2+1\}.
$$
Both values occur, so there are $r$ fibers equal to $\cI$
and $b-r$ equal to $\mathcal J$, where $0<r<b$.
Moreover, $1\le k_1\le d_1$ and $0\le k_2<d_2$.

Let $C_-$ and $C_+$ denote the contributions to
$C(\mathcal A,\tau)$ from negative and positive first
coordinates, respectively.
For each $1\le j<k_1$, there are $b$ sources with first
coordinate $-j$, and at most $r$ of them can map to a target
with first coordinate $-k_1$.
All other targets have first coordinate at least $-k_1+1$.
Thus their total contribution is at least
$$
    r\left\lceil\frac{pj-k_1}{d_1}\right\rceil
    +(b-r)\left\lceil\frac{pj-(k_1-1)}{d_1}\right\rceil.
$$
Adding the contribution of the $r$ sources with first
coordinate $-k_1$ gives
$$
    C_-
    \ge rF_{d_1}(k_1)+(b-r)F_{d_1}(k_1-1).
$$
Similarly, there are only $b-r$ targets with first coordinate
$k_2+1$. For each $1\le j\le k_2$, at least $r$ of the $b$
sources with first coordinate $j$ therefore map to targets
with first coordinate at most $k_2$.
Including the $b-r$ sources with first coordinate $k_2+1$,
we obtain
$$
    C_+
    \ge rF_{d_2}(k_2)+(b-r)F_{d_2}(k_2+1).
$$
The zero-source contributions are nonnegative.
Both $\cI$ and $\mathcal J$ satisfy the hypotheses of
\cref{P:small-one}, so
\begin{align*}
    C(\mathcal A,\tau)
    &\ge r\bigl(F_{d_1}(k_1)+F_{d_2}(k_2)\bigr)
      +(b-r)\bigl(F_{d_1}(k_1-1)+F_{d_2}(k_2+1)\bigr)\\
    &=rS^p(\cI)+(b-r)S^p(\mathcal J)
      =\sum_m S^p(\cI_m)
      \ge bS^p(k).
\end{align*}

\smallskip
\noindent\emph{A strict reduction for a non-weight-minimizing fiber.}
Suppose now that some fiber is not weight-minimizing.
Every weight-minimizing $k$-set has the same maximal weight,
which is a lower bound for the maximal weight of any $k$-set.
Consequently, we can choose a non-weight-minimizing fiber
$\cI_m$ containing an element $x$ whose weight is maximal
among all first coordinates in $\mathcal A$.
Since $\cI_m$ is not weight-minimizing, there is a
$z\notin\cI_m$ with $w(z)<w(x)$.

Replace $(x,m)$ by $(z,m)$ to obtain $\mathcal A'$.
If $(x,m)$ is fixed by $\tau$, replace it by the fixed point
$(z,m)$. Otherwise, remove $(x,m)$ from its cycle, join its
predecessor directly to its successor, and add $(z,m)$ as a
fixed point. Denote the resulting permutation by $\tau'$.
We show that
$C(\mathcal A',\tau')<C(\mathcal A,\tau)$.

The values of $w$ lie in $D^{-1}\Z_{\ge0}$, and $p>2D$
implies $p-1\ge2D$. Hence
\begin{equation}
    \label{E:floor-inequality}
    \lfloor(p-1)w(x)\rfloor-\lceil(p-1)w(z)\rceil
    >(p-1)\bigl(w(x)-w(z)\bigr)-2\ge0.
\end{equation}
If $(x,m)$ was fixed, then
\begin{align*}
    C(\mathcal A,\tau)-C(\mathcal A',\tau')
    &=M((p-1)x)-M((p-1)z)\\
    &=\lceil(p-1)w(x)\rceil-\lceil(p-1)w(z)\rceil>0.
\end{align*}

Otherwise, put
$u=\pi_1(\tau^{-1}(x,m))$ and $v=\pi_1(\tau(x,m))$.
Subadditivity of $M$, which follows by concatenating
representations, gives
$$
    M(pu-v)\le M(pu-x)+M(x-v).
$$
Therefore
\begin{align*}
    C(\mathcal A,\tau)-C(\mathcal A',\tau')
    &=M(px-v)+M(pu-x)-M(pu-v)-M((p-1)z)\\
    &\ge M(px-v)-M(x-v)-M((p-1)z).
\end{align*}
Since $w(v)\le w(x)$ and $x\ne0$, both $px-v$ and $x-v$
have the sign of $x$, allowing $x-v=0$.
Set $d=d_1$ if $x<0$ and $d=d_2$ if $x>0$.
Then
\begin{align*}
    M(px-v)-M(x-v)
    &=\left\lceil\frac{|px-v|}{d}\right\rceil
      -\left\lceil\frac{|x-v|}{d}\right\rceil\\
    &\ge\left\lfloor
      \frac{|px-v|-|x-v|}{d}\right\rfloor
      =\lfloor(p-1)w(x)\rfloor.
\end{align*}
Together with \eqref{E:floor-inequality}, this proves the
strict decrease in the remaining case as well.

Repeat this operation while a non-weight-minimizing fiber
remains. The process terminates because the cost is a
nonnegative integer and strictly decreases at each step.
The resulting pair $(\mathcal A^\circ,\tau^\circ)$ has only
weight-minimizing fibers, so the first part of the proof gives
$$
    C(\mathcal A,\tau)
    >C(\mathcal A^\circ,\tau^\circ)
    \ge bS^p(k).
$$
Thus equality in the asserted lower bound is possible only
when every fiber is weight-minimizing; in that case the first
part also gives the stronger bound by
$\sum_m S^p(\cI_m)$.
\end{proof}

We follow up with some important consequences of this lemma.

\begin{proposition}
    \label{P:cI-is-linspace}
    If $1 \le k \le d_1 + d_2+1$ and $p > 2D$,
    then
    \[
        S^p(k) = S^p(\cI) = 
        \sum_{\ell=1}^{k_1} \left\lceil \frac{p\ell -
            k_1}{d_1}\right\rceil 
            + 
            \sum_{\ell = 1}^{k_2}
            \left\lceil \frac{p\ell-k_2}{d_2} \right\rceil,
    \]
    where $\cI = \{-k_1, \dots, k_2\}$
    is some weight-minimizing set.
\end{proposition}
\begin{proof}
    Let $\cI \in \binom{\Z}{k}$, and let $b = 1$.
    Setting $\mathcal{A} = \cI \times \{0\}$
    and letting $\sigma$ denote any permutation on $\cI$,
    $\sigma$ lifts to an obvious permutation 
    $\tau : \mathcal{A} \to \mathcal{A}$
    by setting $\tau((\ell, 0)) = (\sigma\ell, 0)$.
    \cref{L:key} now implies
    \[
        \sum_{a \in \mathcal{A}} M(p \pi_1(a) - \pi_1(\tau a))
        = \sum_{\ell \in \cI} M(p \ell - \sigma\ell)
        = S^p(k)
    \]
    only if $\cI$ is weight-minimizing,
    otherwise it is greater.
    Therefore if $S^p(k) = S^p(\cI)$,
    then $\cI$ is weight-minimizing.
    In this case, $\cI$ is an interval
    by \cref{L:weight-properties} (2),
    so the formula for $S^p(k) = S^p(\cI)$
    is given by \cref{P:small-one}.
\end{proof}
\begin{corollary}
    \label{C:uniqueness-is-invariant}
    Let $p, p' > 2D$, and let $k \le d_1 + d_2+1$.
    If $p \equiv p' \pmod{d_1}$ and $p \equiv p' \pmod{d_2}$,
    then $S^p(k)$ and $S^{p'}(k)$ are achieved for the same intervals $\cI$.
\end{corollary}
\begin{proof}
    By \cref{P:cI-is-linspace}, every subset achieving $S^p(k)$
    is a weight-minimizing interval $\cI = \{-k_1, \dots, k_2\}$.
    Choose one such subset.
    Because $d_1 \mid p' - p$, we have
    $\left\lceil \frac{p'\ell - k_1}{d_1}\right\rceil = \left\lceil \frac{p\ell
        - k_1}{d_1}\right\rceil + \frac{(p'-p)\ell}{d_1}$,
    and similarly for $d_2$.
    Therefore,
    \[
        S^{p'}(\cI) = S^p(\cI) + (p'-p)W_k.
    \]
    Each candidate minimizing interval has weight $W_k$.
    Its cost therefore changes by the same quantity $(p'-p)W_k$,
    so the minimizing interval or intervals remain unchanged.
\end{proof}

We now establish two results needed for
characterizing the generic polygon defined in 
\cref{S:extension_theorem} and comparing it with the
Newton polygon of the Dwork operator.
\cref{C:kplus1} establishes the convexity of the sequence
$(S^p(k))_{k=1,\dots,d_1+d_2}$,
and identifies strict convexity with the existence of a unique minimizing 
index set.
\cref{P:minimal-bijections} extends the minimization problem from permutations
of a single set to chains of bijections between different sets.
Its lower bound controls the terms in the Cauchy-Binet--type expansion
encountered in \cref{P:generic-np-extension}.

\begin{corollary}[Convexity of $S^p(k)$]
    \label{C:kplus1}
    Assume $1 < k \le d_1 + d_2$ and $p > 2D$.
    Then 
    \[
        2S^p(k) \le S^p(k-1) + S^p(k+1)
    \]
    with equality
    if and only if
    there exist two distinct subsets 
    $\cI_1, \cI_2 \in \binom{\Z}{k}$
    for which $S^p(\cI_1) = S^p(\cI_2) = S^p(k)$.

    Consequently, for every $k \in \{1, \dots, d_1+d_2-1\}$,
    the point $\left(k, \frac{S^p(k)}{p-1}\right)$
    lies on $\GL(d_1,d_2,p)$,
    and $k \in V_p$ (that is, there is an internal vertex
    for $\GL(d_1,d_2,p)$ at $k$)
    if and only if there is a unique $\cI \in \binom{\Z}{k}$
    for which $S^p(k) = S^p(\cI)$.
\end{corollary}
\begin{proof}
    For a positive integer $d < p/2$, let 
    $F_d(n) = \sum_{\ell = 1}^n \left\lceil \frac{p\ell - n}{d}\right\rceil$,
    and define $\Delta_d(j) = F_d(j) - F_d(j-1)$ for $j \ge 1$
    and $\Delta_d(0) = 0$.
    For $j \le d$,
    we have 
    \[
        \Delta_d(j) = \left\lceil \frac{(p-1)j}{d} \right\rceil 
        - \# \left\{1 \le i < j \mid pi \equiv j \pmod{d}\right\},
    \]
    and the last term is always either $0$ or $1$ since $p$ and $d$
    are coprime.  Therefore $\Delta_d(j)$ is positive.
    It also follows that for $j < d$,
    \begin{equation}
        \label{E:Delta-delta}
        \Delta_d(j+1) - \Delta_d(j) \ge 
        \left\lfloor\frac{p-1}{d}\right\rfloor - 1 \ge 1.
    \end{equation}

    Consider the multiset
    \[
        \mathbf{U} \coloneqq 
        \left\{\Delta_{d_1}(j) \mid j=1,2, \dots, d_1 \right\}
        \sqcup 
        \left\{\Delta_{d_2}(j) \mid j=1,2, \dots, d_2 \right\}.
    \]
    Sort the finite elements as $u_1 \le u_2 \le \dots \le u_{d_1+d_2}$.
    Then we have
    by \cref{P:cI-is-linspace}
    that $S^p(k) = u_1 + \dots + u_{k-1}$, since this finite series telescopes.
    Therefore
    \[
        S^p(k+1) - S^p(k) = u_k \ge u_{k-1} = S^p(k) - S^p(k-1),
    \]
    with equality if and only if $u_k = u_{k-1}$.
    Since $p > 2D$,
    the sequences $(\Delta_{d_1}(j))_{j=1}^{d_1}$ and 
    $(\Delta_{d_2}(j))_{j=1}^{d_2}$ are both strictly increasing
    by \eqref{E:Delta-delta}.
    Therefore $u_{k-1} = u_k$
    if and only if for some nonnegative integers $k_1$ and $k_2$
    that sum to $k$,
    $\Delta_{d_1}(k_1) = \Delta_{d_2}(k_2) = u_{k-1} = u_k$.
    This equality is equivalent to the equality
    \[
        \sum_{\ell=1}^{k_1} \left\lceil \frac{p\ell - k_1}{d_1}\right\rceil 
        + \sum_{\ell=1}^{k_2-1} \left\lceil \frac{p\ell -
            (k_2-1)}{d_2}\right\rceil 
        = \sum_{\ell=1}^{k_1-1} \left\lceil \frac{p\ell -
            (k_1-1)}{d_1}\right\rceil + 
            \sum_{\ell=1}^{k_2} \left\lceil \frac{p\ell - k_2}{d_2}\right\rceil.
    \]
    By \cref{P:small-one},
    the lefthand side is equal to $S^p(\cI_1)$,
    where $\cI_1 = \{-k_1, \dots, k_2-1\}$,
    and the righthand side is equal to $S^p(\cI_2)$,
    where $\cI_2 = \{-k_1+1, \dots, k_2\}$,
    and both are equal to $S^p(k) = u_1 + \dots + u_{k-1}
    = u_1 + \dots + u_{k-2} + u_k$.

    The rest of the corollary is immediate
    for all $k>1$,
    and for $k=1$, it's straightforward to observe that $S^p(1)$
    is always a vertex, and its unique minimizing subset is $\{0\}$.
\end{proof}

\begin{proposition}[Cyclic cost bound]
    \label{P:minimal-bijections}
    Let $1 \le k \le d_1 + d_2+1$ and $p > 2D$.
    For two subsets $\cI, \mathcal{J} \in \binom{\Z}{k}$,
    we define
    \[
        S^p(\cI, \mathcal{J}) \coloneqq \min_\sigma 
        \sum_{\ell \in \cI} M(p\ell - \sigma \ell),
    \]
    where $\min_\sigma$ ranges over all bijections 
    $\sigma : \cI \to \mathcal{J}$.

    Let $b \in \Z_{\ge 1}$.
    Then for any
    $(\cI_0, \dots, \cI_{b-1}) \in \binom{\Z}{k}^{\Z / b\Z}$,
    \[
        \sum_{m \in \Z / b\Z} S^p(\cI_{m+1}, \cI_m) \ge bS^p(k).
    \]
    Moreover, if $\cI$ is the unique subset in $\binom{\Z}{k}$ for which
    $S^p(\cI) = S^p(k)$,
    then
    \[
        (\cI_0, \dots, \cI_{b-1}) = (\cI, \dots, \cI)
    \]
    is the unique $b$-tuple for which the above inequality is an equality.
\end{proposition}
\begin{proof}
    Form the disjoint union
    \[
        \mathcal{A} = \bigsqcup_{m \in \Z / b\Z} \cI_m \times \{m\}.
    \]
    For each $m$, let $\sigma_m : \cI_m \to \cI_{m-1}$
    be a minimizing bijection.
    Then we construct a permutation 
    $\tau : \mathcal{A} \to \mathcal{A}$
    given by $\tau(\ell, m) = (\sigma_m(\ell), m-1)$.
    We then have
    \[
        \sum_{m \in \Z / b\Z} S^p(\cI_{m+1}, \cI_m)
        = \sum_{a \in \mathcal{A}} M(p\pi_1(a) - \pi_1(\tau a))
        \ge b S^p(k)
    \]
    by \cref{L:key}, with equality only if $\cI_m$
    is weight-minimizing for all $m$.
    In this case, we also have
    \[
        \sum_{m \in \Z / b\Z} S^p(\cI_{m+1}, \cI_m) 
        \ge \sum_{m \in \Z / b\Z} S^p(\cI_m).
    \]
    If there is only one $\cI \in \binom{\Z}{k}$ that minimizes
    $S^p(\cI)$,
    then of course equality only holds if 
    $\cI_m = \cI$ for all $m$.
    Then setting $\cI_m = \cI$ for all $m$
    implies $S^p(\cI_{m+1}, \cI_m) = S^p(\cI) = S^p(k)$
    by \cref{P:cI-is-linspace},
    so the desired equality holds.
\end{proof}

\subsection{A total ordering on solution vectors}
\label{SS:ordering}
The minimum costs computed above identify the terms 
that will ultimately contribute to the minimum valuation
in the coefficients of the Fredholm determinant
in \cref{S:ck_formula}.
However, there is the potential for term cancellation
in the determinant expansion.
In order to address this, we must identify a term 
in the expansion that cannot cancel.
We will do this by developing a lexicographic-type ordering
on the monomial vectors, and showing that there is a unique term 
that is minimal under this ordering.

Recall that the representations in \cref{D:cmp}
are nonnegative integer vectors in $\Z_{\ge 0}^{I_0}$.

\begin{definition}
    \label{D:total-order}  
    Denote by $\OO$ the following total order on 
    the vectors in $\Z_{\ge 0}^{I_0}$:
    first, for each $i \in I_0$, let $e_i$ denote the standard
    basis vector $e_i = (\delta_{ij})_{j \in I_0}$,
    where $\delta$ is the usual Kronecker delta function.
    We order these basis elements as follows:
    \[
        e_1 <_\OO e_2 <_\OO \cdots <_\OO e_{d_2} <_\OO 
        e_{-1} <_\OO e_{-2} <_\OO \cdots <_\OO e_{-d_1}.
    \]
    Then, given two vectors 
    $\mathbf{s} = (s_i)_{i\in I_0}$
    and $\bt = (t_i)_{i \in I_0}$,
    $\mathbf{s} <_\OO \bt$ if either 
    \begin{enumerate}
        \item $|\mathbf{s}| < |\bt|$, or 
        \item $|\mathbf{s}| = |\bt|$ and 
            $s_i > t_i$ for the smallest $e_i$
            for which $s_i \ne t_i$.
    \end{enumerate}
\end{definition}
We first identify the minimal solution for a particular $N$,
and then minimize the sum of the resulting distinguished vectors over
permutations of an interval.
\begin{lemma}[Minimal solution under $\OO$]
    \label{L:greedy-best}
    For any integer $N$,
    the unique solution in $\Sol(N)$
    of minimal order under 
    $\OO$ is the solution $\bt^+$,
    as defined in \cref{L:cmp-facts}.
\end{lemma}
\begin{proof}
    Let $\bt \in \Sol(N)$,
    and let $|N| = qd+r$,
    where $d = d_1$ if $N<0$ and $d=d_2$
    if $N\ge0$, and where $r \in \{0,\dots,d-1\}$.
    If $N=0$ or $r=0$, then the solution $\bt^+$
    is unique, which trivializes the lemma,
    so suppose $N \ne 0$ and $r>0$,
    and let $\varepsilon = \sgn(N)$.
    Then $|\bt| = q+1$.
    Now let $j \in I_0$.
    If $\sgn(j) \ne \varepsilon$, then $t_j = 0$
    by \cref{L:cmp-facts},
    so suppose $\sgn(j) = \varepsilon$.
    If $|j| < r$,
    then $t_j = 0$,
    because otherwise we would have the contradiction
    $|N| = \sum_{i \in I_0} |i| t_i < qd+r = |N|$.
    If $|j| = r$,
    then $t_j \le 1$ because $2r + (q-1)d < qd+r = |N|$.
    And if $t_j = 1$ then all other coins must have denomination 
    $\pm d$, implying that $\bt = \bt^+$.

    The greedy solution $\bt^+$ has $t^+_{\varepsilon r} = 1$.
    Thus at the earliest possible coordinate $e_{\varepsilon r}$,
    it has a coefficient of $1$,
    while every competing minimal representation has a coefficient of $0$.
    Therefore, $\bt^+$ is uniquely minimal under $\OO$.
\end{proof}

\begin{lemma}[Minimal permutation under $\OO$]
    \label{L:unique-minimal}
    Suppose $-d_1 \le -k_1 \le 0 \le k_2 \le d_2$,
    and let $\cI = \{-k_1, \dots, k_2\}$.
    Furthermore, let $\bt^+(N)$ denote the solution in $\Sol(N)$
    as defined in \cref{L:cmp-facts}.
    If $p$ is coprime to both $d_1$ and $d_2$,
    then there is a unique permutation $\tau$
    for which the summation
    \[
        \sum_{\ell \in \cI} \bt^+(p\ell - \tau\ell)
    \]
    is minimal under $\OO$.
\end{lemma}
\begin{proof}
    We define $\tau$ the same way as in the proof of 
    \cref{P:small-one}.
    Starting with $\ell = k_2$ and working downward inductively,
    we let $I_\ell = \{1, \dots, k_2\} \setminus\{\tau^{-1}(j) \mid j \ge
    \ell+1\}$.
    Then we set $\tau^{-1}(\ell) = i$,
    where $i$ denotes the element in $I_\ell$ that minimizes $[pi-\ell]_{d_2}$.
    For negative $\ell$,
    the same construction is done with $(k_1, d_1)$ instead of $(k_2, d_2)$,
    and with flipped signs.
    Specifically, let $\tau'$ be a permutation of $\{1, \dots, k_1\}$,
    determined iteratively by setting $\tau'^{-1}(-\ell) = i$,
    where $i$ denotes the element in $\{1, \dots, k_1\}
    \setminus\{\tau'^{-1}(j)\mid j \ge -\ell+1\}$
    that minimizes $[pi+\ell]_{d_1}$.
    Then set $\tau(\ell) = -\tau'(-\ell)$,
    and fix $\tau(0) = 0$.
    
    Note that in the proof of \cref{P:small-one},
    we already established that
    $\tau^{-1}$ is well-defined and injective,
    and that $\tau \in \bSig^p(\cI)$.

    Let $\sigma$ be another permutation,
    and suppose
    \begin{equation}
        \label{E:lex-inequality}
        \sum_{\ell \in \cI} \bt^+(p\ell - \sigma\ell)
        \le_{\OO} \sum_{\ell \in \cI} \bt^+(p\ell - \tau\ell).
    \end{equation}
    We will show that $\sigma = \tau$ on the positive elements by a reverse
    induction argument on $\sigma^{-1}$ and $\tau^{-1}$.
    Let $\ell \in \{1, \dots, k_2\}$, and assume $\sigma^{-1}(j) = \tau^{-1}(j)$
    for all $j > \ell$.

    Let $z = \tau^{-1}(\ell)$, and let $r = [pz - \ell]_{d_2}$.
    If $r = 0$, then we must have $\sigma^{-1}(\ell) = z$,
    as otherwise, we will have $\sigma z < \ell$, implying
    $\left\lceil \frac{pz - \sigma z}{d_2}\right\rceil > 
    \left\lceil \frac{pz - k_2}{d_2}\right\rceil$,
    and this will contradict \cref{P:small-one}.
    So assume $r > 0$.

    Let $y \ge 1$, and suppose $\sigma y \le \ell$.
    If $[py-\sigma y]_{d_2} = 0$, then let $r' = d_2$;
    otherwise let $r' = [py - \sigma y]_{d_2}$.
    We claim that either $r' > r$ or $\sigma^{-1}(\ell) = \tau^{-1}(\ell)$.

    Note that if $y = \tau^{-1}(j)$ for some $j \ge \ell+1$,
    then our inductive hypothesis establishes $y = \sigma^{-1}(j)$,
    which runs contrary to our assumption that $\sigma y \le \ell$.
    Thus $y \in I_\ell$, and by the minimality of $z = \tau^{-1}(\ell)$,
    we must have that
    \begin{equation}
        \label{E:py-inequality}
        [py - \ell]_{d_2} \ge [pz - \ell]_{d_2} = r
    \end{equation}
    with equality if and only if $y = z$.
    On the other hand, 
    since $\OO$ first compares total degree
    and $\tau \in \bSig^p(\cI)$,
    it follows from
    \eqref{E:lex-inequality}
    that $\sigma \in \bSig^p(\cI)$,
    thus by \cref{P:small-one},
    $\left\lceil \frac{py - \sigma y}{d_2}\right\rceil =
    \left\lceil \frac{py - k_2}{d_2}\right\rceil$.
    Thus $[py-j]_{d_2} \ne 0$ for any $j \in \{\sigma y + 1, \dots, k_2\}$,
    and is therefore decreasing on this interval.
    Thus either $[py-\sigma y]_{d_2} = 0$,
    or $[py - \sigma y]_{d_2} \ge [py - \ell]_{d_2}$,
    with equality if and only if $\sigma y = \ell$.
    Combining this with \eqref{E:py-inequality},
    we obtain $r' \ge r$, with equality if and only if $y = z$
    and $\sigma y = \ell$.
    Since $y=z$ and $\sigma y = \ell$ is equivalent to $\sigma^{-1}(\ell) =
    \tau^{-1}(\ell)$, our claim is proven.

    If $r' \ne d_2$, then the coefficient of $e_{r'}$ in $\bt^+(py - \sigma y)$
    is equal to $1$.
    All other coefficients, aside from the coefficient of $e_{d_2}$, are zero.
    If $r' = d_2$, then all coefficients of $\bt^+(py-\sigma y)$ except possibly that of $e_{d_2}$ vanish.
    Therefore, for any $1 \le s \le r$,
    the coefficient of $e_s$ in 
    $\sum_{\{y>0 \mid \sigma y \le \ell\}} \bt^+(py - \sigma y)$ 
    is equal to $0$,
    unless $s = r$ and $\sigma^{-1}(\ell) = \tau^{-1}(\ell)$.
    By \eqref{E:lex-inequality},
    this must be the case.

    A symmetric argument shows that $\sigma = \tau$
    on the negative elements of $\cI$ as well,
    and we conclude $\tau$ to be minimal under $\OO$.
\end{proof}
 \section{Dwork theory and exponential-sum 
\texorpdfstring{$L$}{L}-functions}
\label{S:dwork-work}
In this section, we lay the groundwork for future sections by
defining the Dwork operator with respect to 
a Laurent polynomial $f$.
The Dwork operator is a completely continuous
endomorphism whose Fredholm determinant possesses a Newton polygon that matches
the normalized Newton polygon
of the $L$-function of $f$ on the interval $[0,d_1+d_2]$.
We will express the Dwork operator as a matrix indexed by 
the monomial formal basis $\{s^i \mid i \in \Z\}$, 
and prove some properties about its
entries with the language of change-making problems.

\subsection{The Dwork operator}
\label{SS:dwork-operator}
The power series
\begin{equation}
    \label{E:AHE}
    E_p(s) \coloneqq \exp \left( \sum_{m=0}^\infty \frac{s^{p^m}}{p^m} \right)
    = \sum_{i=0}^\infty \beta_i s^i
\end{equation}
is the $p$-adic Artin-Hasse exponential function.
It lies in $(\Q \cap \Z_p)[[s]]$.
Fix a root $\gamma$
of $\log(E_p(s))$
with $p$-adic valuation $v_p(\gamma) = \frac{1}{p-1}$.
$E_p(\gamma)$ is a primitive $p$-th root of unity,
which we shall denote as
$\zeta_p$.
We have that $\Z_p[\zeta_p] = \Z_p[\gamma]$.

Let $f = f(x) = \sum_{i=-d_1}^{d_2} a_i x^i \in \F_q[x^{\pm1}]$.
 Assume $a_{-d_1}a_{d_2} \ne 0$,
so that $f$ has fixed positive degree $d_2$ and negative degree $d_1$.
Also assume that $p \nmid d_1d_2$.
For each $m \ge 1$, the $m$-th exponential sum $S_m^*(f)$
is given by
\begin{equation}
    \label{E:exponential-sum}
    S_m^*(f) \coloneqq \sum_{z \in \F_{q^m}^\times}
    \zeta_p^{\Tr_{\F_{q^m}/\F_p}(f(z))},
\end{equation}
and the $L$-function of the exponential sum of $f$
is given by 
\begin{equation}
    \label{E:L}
    L^*(f; s) = \exp \left(\sum_{m=1}^\infty S^*_m(f) \frac{s^m}{m}\right).
\end{equation}
It is a well known corollary of the Weil Conjecture for curves that 
the $L$-function is a polynomial in $1 + s\Z_p[\zeta_p][s]$,
of degree $d_1 + d_2$;
see \cite[Theorem 1.14]{AS87}.
Moreover,
the valuations of the reciprocal roots of $L^*(f; s)$
are symmetric, so that for every line segment $S$ with length $\ell$
and slope $\lambda$ in
the $q$-adic Newton polygon of $L^*(f; s)$,
there is a corresponding line segment $S'$ with length $\ell$
and slope $1 - \lambda$.
Since $L^* \in 1 + s \Z_p[\zeta_p][s]$, it follows that 
every such $\lambda$ must lie between $0$ and $1$.

\begin{remark}
    \label{r:assume-a0-zero}
    Let $f^* = f - a_0$.
    If we set $z_0 = -\Tr_{\F_q / \F_p}(a_0)$,
    then $S_m^*(f^*) = \zeta_p^{m z_0} S_m^*(f)$,
    and $L^*(f^*; s) = L^*(f; \zeta_p^{z_0}s)$.
    Thus the Newton polygon of $L^*(f^*; s)$
    agrees with the Newton polygon of $L^*(f; s)$,
    and we can therefore assume that $a_0 = 0$
    for our purposes.
\end{remark}

To the Laurent polynomial $f$ we attach the power series
\begin{equation}
    \label{E:Ft}
    F(s) = \prod_{i\in I_0} E_p(\wh{a}_i \gamma s^i)
    = \sum_{n=-\infty}^\infty G_n s^n,
\end{equation}
where $\wh{a}$ denotes the Teichm\"uller representative of $a \in \F_q$
in $\Q_q$, the unique degree $b$ unramified extension of $\Q_p$.

Let $w : \Z \to \Q_{\ge 0}$ be given by
\[
    w(i) = 
    \begin{cases*}
        \frac{-i}{d_1} & if $i \le 0$ \\
        \frac{i}{d_2} & if $i > 0$.
    \end{cases*}
\]
Then
$F(s)$ lies in the space
\[
    \mathcal{H} \coloneqq \left\{\sum_{j=-\infty}^\infty B_j s^j 
        \in \Q_q(\zeta_p)[[s^{\pm 1}]]
        \:\Big|\: 
        v_p(B_j) \ge \frac{w(j)}{p-1} + O(1)\right\}.
\]
When equipped with the norm 
\[
    \left|
    \sum_j B_j s^j
    \right| = 
    \sup_{j \in \Z}\left(|B_j|_p p^{w(j)/(p-1)}\right),
\]
$\mathcal{H}$ is the $\Q_q(\zeta_p)$-Banach algebra
consisting of power series $\xi \in \Q_q(\zeta_p)[[s^{\pm1}]]$
for which $|\xi| < \infty$.

Let $\tau$ denote the lifting of the Frobenius automorphism on $\F_q$
to $\Q_q$.
The canonical Dwork operator,
as in \cite{ref:bombieri},
is given first by setting
$    \alpha_0 = \tau^{-1} \circ \psi_p \circ F,
$
where $\psi_p \left(\sum_i B_i s^i\right) = \sum_i B_{pi} s^i$
and $F(H) = F \cdot H$,
and then by setting
\begin{equation}
    \label{E:alpha}
    \alpha = \alpha_0^b.
\end{equation}
Following the theory of Serre in 
\cite{ref:serre1965},
the operator $\alpha$
is completely continuous on the Banach space $\mathcal{H}$.

We omit the proof of the following proposition 
as it is well known, see \cite{ref:bombieri}.

\begin{proposition}
    \label{prop:L-fredholm-NP}
    Let $f$ be a Laurent polynomial defined over $\F_q$,
     with positive degree $d_2$ and negative degree $d_1$,
    and $p\nmid d_1d_2$.
    Let  $L^*(f; s)$ and $\alpha$
    be as in \eqref{E:L} and \eqref{E:alpha},
    respectively.
    Then
    \[
        \NP_q(L^*(f; s)) = \NP_q(\det(1 - \alpha s))\big|_{[0,d_1+d_2]}.
    \]
\end{proposition}
There exists a matrix expression for $\alpha$
under the formal basis $\cB = s^{\Z} = \{s^i \mid i \in \Z\}$.
Let $\cG = \Mat_{\cB}(\psi_p \circ F)$ denote
the matrix for $\psi_p \circ F$ indexed by $\cB$.
Then one easily verifies that
\begin{equation}
    \label{E:matG-index}
    \cG = (G_{pi-j})_{i,j \in \Z}.
\end{equation}
Write
\begin{equation}
    \label{E:matGb}
    [\cG, \tau, b] = \cG^{\tau^{b-1}} \cG^{\tau^{b-2}} \cdots \cG^\tau \cG.
\end{equation}
Let $\Mat_{\cB}(\alpha_0)$ denote the matrix obtained by
setting the $i$-th column equal to $\alpha_0(s^i)$.
Then $\cG^{\tau^{-1}} = \Mat_{\cB}(\alpha_0)$,
and since $\tau^b = 1$,
it follows that 
\begin{equation}
    \label{E:mat-alpha}
    \Mat_{\cB}(\alpha) = 
    \Mat_{\cB}(\alpha_0^b) = \cG^{\tau^{-1}} \cG^{\tau^{-2}} \cdots 
    \cG^{\tau^{-b}} = [\cG, \tau, b].
\end{equation}
With \eqref{E:mat-alpha},
it follows from \cref{prop:L-fredholm-NP}
that
\begin{equation}
    \label{E:L*-G-tau-b}
    \NP_q(L^*(f; s)) = \NP_q(\det(1 - [\cG, \tau, b]s))|_{[0,d_1+d_2]}.
\end{equation}
In a later section (see \cref{S:extension_theorem}),
we will prove that 
the valuations of the interior vertices of $\det(1 - \cG s)$
will characterize when $\det(1 - [\cG, \tau, b]s)$
attains the generic Newton polygon.
Thus we finish with a statement about the former
Fredholm determinant.
\begin{proposition}
    \label{P:cG-det-expression}
    Then $\det(1 - \cG s) = 1 + C_1 s + C_2 s^2 + \dots$,
    where 
    \[
        C_k = (-1)^k \sum_{(\cI, \sigma)} \sgn(\sigma) 
        \prod_{\ell \in \cI} G_{p\ell - \sigma \ell},
    \]
    where the sum $\sum_{(\cI, \sigma)}$ ranges over
    all pairs $(\cI, \sigma)$ of subsets
    $\cI \in \binom{\Z}{k}$ 
    and permutations $\sigma$ of $\cI$.
\end{proposition}

\subsection{Rank-\texorpdfstring{$\ell$}{l} extensions}
\label{SS:rank-l}
We introduce
the $L$-functions associated with additive characters of $\F_{p^\ell}$
and reduce them to the rank-one case.

The function $\chi_1 : \F_p \to \overline{\Q}_p^\times$ given by $\chi_1(z) =
\zeta_p^{z}$
is an additive character of $\F_p$.
Equations \eqref{E:exponential-sum} and \eqref{E:L}
correspond to the exponential sums and $L$-functions associated with this
particular nontrivial additive character.
If $\ell \in \Z_{\ge 1}$ is a factor of $b$, then we may choose a nontrivial
additive character $\chi_\ell : \F_{p^\ell} \to \overline{\Q}_p^\times$.
For this character, we have analogues 
of \eqref{E:exponential-sum} and \eqref{E:L}
by defining the $m$-th exponential sum for $f$ over $\F_q$ corresponding to
$\chi_\ell$:
\[
    S_m^*(f, \chi_\ell) = \sum_{z \in \F_{q^m}^\times} 
    \chi_\ell(\Tr_{\F_{q^m} / \F_{p^\ell}}(f(z)))
\]
and its $L$-function
\[
    L^*(f, \chi_\ell; s) = 
    \exp \left(\sum_{m=1}^\infty S_m^*(f, \chi_\ell) \frac{s^m}{m}\right).
\]
This is the rank-$\ell$ $L$-function of $f$,
and it can be reduced to a rank-$1$ $L$-function, as in the following
proposition.
\begin{proposition}
    \label{P:rank-l-reduction}
    Let $\ell \mid b$.
    For every nontrivial additive character 
    $\chi_\ell : \F_{p^\ell} \to \overline{\Q}_p^\times$,
    there exists a constant $\alpha \in \F_{p^\ell}^\times$
    for which, for all Laurent polynomials 
    $f \in \F_q[x^{\pm 1}]$, we have 
    \[
        L^*(f, \chi_\ell; s) = L^*(\alpha f; s).
    \]
\end{proposition}
See the discussion above \cite[Proposition 2.3]{zhu2025}
for a proof of the statement for $f \in \F_q[x]$ using the 
trace-dual basis; 
the same argument can be applied for this case.

On the other hand, the rank-$\ell$ $L$-function relates
to the zeta function of $X_{f, \ell} : y^{p^\ell} - y = f$
in the following well-known manner,
\begin{proposition}
    \label{P:zeta-L-connection}
   Let $\ell \mid b$.
    Let $f \in \F_q[x^{\pm 1}]$ be a Laurent polynomial
    with poles at $0$ and $\infty$,
    each with order coprime to $p$.
    If $X_{f, \ell} : y^{p^\ell} - y = f(x)$ denotes 
    the rank-$\ell$ Artin--Schreier curve associated to $f$,
    then
    \[
        Z(X_{f, \ell} / \F_q; s) = 
        \frac{\prod_{\lambda \in \F_{p^\ell}^\times}
        L^*(\lambda f; s)}{(1-s)(1-qs)}.
    \]
\end{proposition}

\subsection{\texorpdfstring{$T$}{T}-adic exponential sums and 
\texorpdfstring{$C$}{C}-functions}
\label{SS:T-adic-intro}

Liu and Wan developed $T$-adic exponential sums to interpolate exponential sums associated with characters of $p$-power order. Here we recall the definitions needed below. 
For $g \in \F_q[x^{\pm 1}]$,
write 
\[
    S_m^*(g ; T) \coloneqq \sum_{x \in \F_{q^m}^\times} 
    (1 + T)^{\Tr_{\Q_{q^m} / \Q_p}(\wh{g}(\wh{x}))}
\]
and 
\[
    C^*(g; T, s) \coloneqq \exp \left(-\sum_{m = 1}^\infty 
    \frac{S_m^*(g ; T)}{q^m-1} \frac{s^m}{m} \right).
\]
Here, $\Q_{q^m} / \Q_p$ is the unramified extension with residue field
$\F_{q^m}$, and 
$\wh{g}$ and $\wh{x}$ denote Teichm\"uller lifts,
with lifts in $\wh{g}$ occurring coefficientwise.
We also define 
\[
    L^*(g; T, s) \coloneqq \frac{C^*(g; T, s)}{C^*(g; T, qs)}.
\]

For each $k \ge 1$,
let $\zeta_{p^k}$ denote a fixed primitive $p^k$-th root of unity,
selecting $\zeta_{p^1}$ so that it agrees with
the primitive root $\zeta_p = E_p(\gamma)$ at the start of this section,
and set $\pi_k = \zeta_{p^k}-1$.
Define orders normalized at $T$, $\pi_k$, and $q$:
\[
    \ord_T(T) = 1, \quad \ord_{\pi_k}(\pi_k) = 1, \quad \ord_q(q) = 1, \quad \ord_{\pi_m}(q)=bp^{m-1}(p-1).
\]
The $T$-adic and $\pi_k$-adic Newton polygons are defined analogously to the
$p$-adic and $q$-adic Newton polygons, with respect to the above valuations.
The $T$-adic valuation is taken in characteristic zero, so that
$\ord_T(p)=0$.
Note that when $k = 1$, we have
\begin{equation}
    \label{E:L-Tadic}
    L^*(g; T, s)|_{T=\pi_1} = L^*(g; s).
\end{equation}

 \section{A formula for the Fredholm coefficients and their valuation}
\label{S:ck_formula}
In \cref{SS:dwork-operator},
we introduced the splitting function $F(s)$,
the matrix $\cG$ associated with the Dwork operator, 
and the Fredholm determinant
$\det(1 - \cG s)$,
with defined expressions
\[
    F(s) = \sum_{n \in \Z} G_n s^n,
    \quad 
    \cG = (G_{pi-j})_{i, j \in \Z},
    \mbox{ and } 
    \det(1 - \cG s)= 1 + \sum_{k=1}^\infty C_k s^k.
\]
(See equation \eqref{E:Ft}, \eqref{E:matG-index}, and 
\cref{P:cG-det-expression}, respectively.)
These coefficients $G_n$ and $C_k$ lie in the ring of integers
of $\Q_q(\zeta_p)$.

In this section, we replace the specific choice of coefficients
$a_i \in \F_q$
with free variables $A_i$.
The analogue of the coefficients $G_n$ of $F(s)$
will be denoted by $\bG_n$,
while the analogue of $C_k$ will be $\bC_k$.
Both $\bG_n$ and $\bC_k$ will then be power series in a Tate algebra
with coefficients in $(\Z_p \cap \Q)[\gamma]$,
where $\gamma$ denotes a root of the logarithm of the $p$-adic Artin-Hasse 
exponential function $E_p(s)$ with valuation $(p-1)^{-1}$.
We can then assign the Gauss norm to this algebra,
inherited from the $p$-adic valuation on the base field.
We will develop the structure of these coefficients in this section,
and express a formula for their valuations
in the language developed in \cref{S:appendix}.

\subsection{Expressions in the coefficient space}
\label{SS:coef-space}
Inherit all notation from \cref{SS:dwork-operator}.
In this first subsection,
we develop the domain space
by replacing the coefficients of $f$ with variables.
We fix the bounds $d_1$ and $d_2$ 
on the negative and positive exponents of $f$, respectively.
Assume both $d_1$ and $d_2$ are positive integers
coprime to $p$.

Define free variables $A_i$ for $i \in I_0$,
and let $A = (A_i)_{i \in I_0}$.
We define the power series
\begin{equation}
    \label{E:FAt}
    F_A(s) = \prod_{i \in I_0} E_p(A_i \gamma s^i)
    = \sum_{n=-\infty}^\infty \bG_n s^n.
\end{equation}
It follows from \cref{P:bG-expression} below that
each $\bG_n = \bG_n(A)$ is a power series in
the Tate algebra $\Z_p[\gamma]\left<A\right>$
equipped with the standard Gauss norm,
and that all coefficients are in $\Q[\gamma]$.
For any particular $f = \sum_{i \in I_0}a_i x^i$
with coefficients in $\F_q$,
we have
\[
    \bG_n(\wh{a}) = G_n
\]
from \cref{SS:dwork-operator}, where $\wh{a} = (\wh{a}_i)_{i \in I_0}$.

To describe $\bG_n$ structurally, 
we recall \cref{D:cmp} from \cref{SS:cmp}.
A \textit{representation} of an integer $n$
in our context is a tuple 
$\bt = (t_i)_{i \in I_0} \in \Z_{\ge 0}^{I_0}$
for which $n = \sum_{i \in I_0} t_i i$.
The value $M(n)$ denotes the smallest possible
value of $|\bt| = \sum_i t_i$ among all representations of $n$,
and a representation of $n$ that achieves the minimum value is called
a \textit{solution} for $n$.

Let $\cT(n)$ denote the set of all representations of $n$,
and let $\Sol(n)$ denote the set of all solutions for $n$.
\begin{proposition}
    \label{P:bG-expression}
    Let $E_p(s) = \sum_{m=0}^\infty \beta_m s^m$.
    For every $n \in \Z$, we have
    \[
        \bG_n = \sum_{\bt \in \cT(n)} 
        \left(\prod_{i \in I_0} \beta_{t_i}\right)
        \gamma^{|\bt|} A^\bt,
    \]
    where $A^\bt$ denotes the monomial 
    $\prod_{i\in I_0} A_i^{t_i}$.
    Moreover, we have $v_p(\bG_n) \ge \frac{M(n)}{p-1}$,
    and if $-pd_1 < n < pd_2$,
    then $v_p(\bG_n) = \frac{M(n)}{p-1}$.
\end{proposition}
\begin{proof}
    The expression for $\bG_n$ is obtained by collecting the coefficients
    in \eqref{E:FAt}.
    The inequality $v_p(\bG_n) \ge \frac{M(n)}{p-1}$
    is an immediate corollary to this expression,
    since $\prod_{i \in I_0} \beta_{t_i}$ is a $p$-adic integer.

    Suppose $-pd_1 < n < pd_2$.
    Assume first that $n > 0$.
    Consider the solution $\bt^+ \in \Sol(n)$,
    as given in \cref{L:cmp-facts}.
    Since $n < pd_2$,
    $t^+_{d_2} = \left\lfloor n/d_2\right\rfloor < p$,
    and thus $t^+_j < p$ for all $j \in I_0$.
    It follows that $\prod_{i \in I_0} \beta_{t^+_i}$
    is a $p$-adic unit,
    as $\beta_j = \frac{1}{j!} \in \Z_p^\times$
    for all $j < p$.
    Therefore the coefficient of $A^{\bt^+}$ in $\bG_n$
    has valuation $\frac{M(n)}{p-1}$,
    and so $v_p(\bG_n) = \frac{M(n)}{p-1}$.

    If $n < 0$,
    the proof is similar.
    The components of the solution $\bt^+$ will all be less than $p$,
    and thus the coefficient of $A^{\bt^+}$ divided by $\gamma^{M(n)}$
will be a $p$-adic unit,
    and so $v_p(\bG_n) = \frac{M(n)}{p-1}$.
    Finally, if $n=0$, then the constant coefficient 
    of $\bG_0$ is equal to $1$, so the equality is immediate.
\end{proof}

Next we extend
the coefficients $C_k$ in \cref{P:cG-det-expression}
into our coefficient space.
For every integer $k \ge 1$, we define
\begin{equation}
    \label{E:bC_k}
    \bC_k(A) \coloneqq (-1)^k \sum_{(\cI, \sigma)} \sgn(\sigma)
    \prod_{\ell \in \cI} \bG_{p\ell - \sigma\ell}(A),
\end{equation}
where $\sum_{(\cI, \sigma)}$ ranges over all pairs
$(\cI, \sigma)$ of sets $\cI \in \binom{\Z}{k}$
and permutations $\sigma$ of $\cI$.
We once again have that $\bC_k$
lies in the Tate algebra $\Z_p[\gamma]\left<A\right>$,
and its coefficients lie in $\Q[\gamma]$.
Moreover, for a particular $f = \sum_{i \in I_0} a_i x^i$,
\[
    \bC_k(\wh{a}) = C_k,
\]
where $C_k$ is the coefficient in \cref{P:cG-det-expression},
and $\wh{a} = (\wh{a}_i)_{i \in I_0}$.

Our goal is to establish a precise value for
$v_p(\bC_k)$, which we will do in 
\cref{C:Ck-circ-nonzero}.

\subsection{Isolating terms of minimal valuation}
\label{SS:Ck-Gn-lowest-terms}
We next isolate the terms in $\bC_k$ with lowest valuation.
Note that for every $n \in \Z$, by \cref{P:bG-expression}
we may express $\bG_n$ as a formal power series over $\gamma$.
Namely,
\[
    \bG_n = \sum_{v = 0}^\infty \bG_{n, v} \gamma^v,
\]
where 
\[
    \bG_{n, v} = \sum_{
        \substack{\bt \in \cT(n)\\ |\bt|=v}}
        \left(\prod_{i\in I_0} \beta_{t_i}\right) A^{\bt}.
\]
Note that $\bG_{n, v} = 0$ for all $v < M(n)$.
Let
\begin{equation}
    \label{E:bGn-min}
    \bG_n^\circ \coloneqq \bG_{n, M(n)}
    = \sum_{\bt \in \Sol(n)} 
    \left(\prod_{i\in I_0} \beta_{t_i}\right) A^{\bt}.
\end{equation}
Note that $\Sol(n)$ is a finite set,
since there are only finitely many tuples $(t_i)_{i\in I_0}$
of nonnegative integers that sum to a fixed integer.
Therefore, $\bG_n^\circ$ is a homogeneous polynomial of 
total degree $M(n)$ in $(\Z_p \cap \Q)[A]$.
It is also nonzero, since every $\beta_{t_i}$ is nonzero.

For each $k \ge 1$, 
we may similarly express $\bC_k$ as a power series over $\gamma$:
\begin{equation}
    \label{E:bCk-power-series}
    \bC_k = \sum_{v=0}^\infty \bC_{k, v} \gamma^v,
\end{equation}
where $\bC_{k, v}$ denotes the collection of monomials in $\bC_k$
with total degree $v$, after factoring out the common term $\gamma^v$.

Here we recall the notation from \cref{D:cycle-sums}
--- namely $S^p(k)$.
This is equal to the minimum possible value of
$\sum_{\ell \in \cI} M(p\ell - \sigma\ell)$,
as $\cI$ and $\sigma$ vary over the same parameter space as in 
\eqref{E:bC_k}.
Let $\bSig_k^p = \left\{(\cI, \sigma) \mid \sum_{\ell \in \cI} M(p\ell -
\sigma\ell) = S^p(k) \right\}$.
Then we may define, analogous to $\bG_n^\circ$,
the term
\begin{equation}
    \label{E:bCk-min}
    \bC_k^\circ \coloneqq 
    (-1)^k \sum_{(\cI, \sigma) \in \bSig_k^p}
    \sgn(\sigma)
    \prod_{\ell \in \cI} \bG_{p\ell - \sigma\ell}^\circ.
\end{equation}
By unraveling the expressions of $\bC_{k,v}$ and $\bG_n$,
one finds
that $\bC_{k, v} = 0$ for all $v < S^p(k)$,
that $\bC_{k, S^p(k)} = \bC_k^\circ$,
and that $\bC_k^\circ$ is a homogeneous polynomial of total degree $S^p(k)$
in $(\Z_p \cap \Q)[A]$.

We summarize this observation with a proposition.

\begin{proposition}
    \label{P:circ-minimal}
    Suppose $k \ge 1$
    and $p > 2D $.
    Then $v_p(\bC_k) \ge (p-1)^{-1} S^p(k)$.
    If $\bC_k^\circ$ is a nonzero 
    polynomial
    modulo $p$,
    then $v_p(\bC_k) = (p-1)^{-1} S^p(k)$.
    \qed
\end{proposition}
We can therefore establish a formula for $v_p(\bC_k)$ 
by establishing that $\bC_k^{\circ}$ is nonzero modulo $p$.
In order to accomplish this,
we will rely heavily on many results from \cref{S:appendix}
to isolate a nonzero term.

\subsection{A nonzero term in 
\texorpdfstring{$\bC_k^\circ$}{Ck°}}
\label{SS:minimal-Ck-nonzero}
Assume $p > 2D$ throughout this subsection.
Recall the combinatorial cost polygon $\GL(d_1,d_2,p)$
from \cref{D:ANP}, along with the set $V_p$
of abscissae of internal vertices of the polygon.
By \cref{C:kplus1},
an internal point $k$ lies in $V_p$
if and only if there is a unique $\cI \in \binom{\Z}{k}$
for which $S^p(k) = S^p(\cI)$.

When $k \in V_p$, therefore,
the formula for $\bC_k^\circ$
loses one layer of complexity.
Because only one minimizing set exists,
\eqref{E:bCk-min} can be expressed as 
\begin{equation}
    \label{E:bCk-min-only}
    \bC_k^\circ = (-1)^k \sum_{\sigma \in \bSig^p(\cI)}
    \sgn(\sigma) \prod_{\ell \in \cI} \bG_{p\ell - \sigma\ell}^\circ,
\end{equation}
where $\bSig^p(\cI) = \{\sigma \in \Sym(\cI) \mid \sum_\ell M(p\ell -
\sigma\ell) = S^p(\cI) = S^p(k)\}$.

Recall the total ordering $\OO$ on the set $\Z_{\ge 0}^{I_0}$
from \cref{SS:ordering}.
This ordering compares the sum of the coordinates first, 
then follows a reverse lexicographic order 
with coordinates ordered (from smallest to largest)
$1, \dots, d_2$, then $-1, \dots, -d_1$.
See \cref{D:total-order} for the complete formal definition.
Extend $\OO$ to the set of monomials in $\Z_p\left<A\right>$
in the natural way.
\begin{proposition}
    \label{P:lexicography}
    If $p>2D$ and $k \in V_p$,
    then $\bC_k^\circ \not \equiv 0 \pmod{p}$.
\end{proposition}
\begin{proof}
    Since $k \in V_p$,
    there is a unique minimizing set $\cI \in \binom{\Z}{k}$,
    and thus we can express $\bC_k^\circ$
    as in \eqref{E:bCk-min-only}.
    Let
    $\Sol(\cI, \sigma) = \prod_{\ell \in \cI} \Sol(p\ell - \sigma\ell)$.
    Then by \cref{L:greedy-best},
    \[
        \min_{(\bt_\ell) \in \Sol(\cI, \sigma)}
        \sum_{\ell \in \cI} \bt_\ell = \sum_{\ell \in \cI} \bt_\ell^+,
    \]
    where $\bt_\ell^+$ denotes the solution to 
    $p\ell - \sigma\ell$ given by \cref{L:cmp-facts},
    and the minimum is under the ordering $\OO$.
    Moreover, the minimum is uniquely achieved at the tuple
    $(\bt_\ell^+)_{\ell \in \cI}$.

    By \cref{P:cI-is-linspace},
    $\cI$ is of the form $\{-k_1, \dots, k_2\}$
    for some $-d_1 \le -k_1 \le 0 \le k_2 \le d_2$,
    and by \cref{L:unique-minimal},
    there is a unique permutation $\tau$ 
    for which
    \[
        \mathbf{T}_{\OO} \coloneqq \sum_{\ell \in \cI} \bt^+(p\ell - \tau\ell)
    \]
    is minimal under $\OO$.
    Its minimality implies $\tau \in \bSig^p(\cI)$,
    and its uniqueness implies that the
    coefficient of $A^{\mathbf{T}_{\OO}}$ in $\mathbf{C}_k^\circ$
    is given by
    \begin{equation}
        \label{E:min-coef}
        D_{\OO} \coloneqq 
        (-1)^k\sgn(\tau) \prod_{\ell \in \cI}
        \prod_{i\in I_0} \beta_{(\bt^+(p\ell - \tau\ell))_i}.
    \end{equation}
    For ease of notation,
    fix $\ell \in \cI$,
    and let $\bt^+(p\ell - \tau\ell) = (t_i)_i$.
    If $i = -d_1$ and $\ell < 0$,
    then $t_i = \left\lfloor \frac{p\ell - \tau\ell}{-d_1}\right\rfloor 
    \le \left\lfloor \frac{(p-1)d_1}{d_1}\right\rfloor = p-1$,
    so $\beta_{t_i}$ is a $p$-adic unit.
    (Recall that $\beta_m$ denotes the coefficient of $x^m$
    in the Artin-Hasse exponential $E_p(x)$).
    Similarly, we find that $t_i \le p-1$
    if $i=d_2$ and $\ell > 0$.
    In all other cases, $t_i = 0$ or $t_i = 1$,
    and $\beta_0 = \beta_1 = 1$.
    Thus, $D_\OO$ is a $p$-adic unit,
    and our proof is complete.
\end{proof}

\begin{corollary}
    \label{C:Ck-circ-nonzero}
    Suppose $p > 2D$ and $k \in V_p$.
    Then 
    for any tuple $z = (z_i)_{i \in I_0}$ of integral elements
    in $\overline{\Q}_p$, we have
    \[
        v_p(\bC_k(z)) = \frac{S^p(k)}{p-1} \iff 
        \bC_k^\circ(z) \not\equiv 0 \pmod{\km_{\overline{\Q}_p}},
    \]
    where $\km_{\overline{\Q}_p}$ denotes the unique maximal ideal in
    ${\mathcal O}_{\overline{\Q}_p}$.
\end{corollary}
\begin{proof}
    By \cref{P:lexicography},
    $\bC_k^\circ \not\equiv 0 \pmod{p}$.
    This establishes $\bC_k^\circ\ne 0$ .
    By
    \cref{P:circ-minimal},
    it now follows that $v_p(\bC_k) = (p-1)^{-1}S^p(k)$,
        The equivalence follows immediately from the definition
    of $\bC_k^\circ$.
\end{proof}

\subsection{A \texorpdfstring{$p$}{p}-independent factorization lemma}
\label{SS:pind-factor}
When a particular Laurent polynomial
$f = \sum_{i \in I_0} a_i x^i \in \F_q[x^{\pm1}]$
is chosen, and the Teichm\"uller lifts
$\wh{a} = (\wh{a}_i)_{i \in I_0}$
of the coefficients are substituted for $A$,
the expression $\bC_k(\wh{a})$ gives the coefficient $C_k$ of $s^k$
in the Fredholm determinant $\det(1 - \cG s)$
associated with $f$, where
\[
    \cG = (G_{pi-j})_{i,j \in \Z} = (\bG_{pi-j}(\wh{a}))_{i, j \in \Z}.
\]
Therefore, by \cref{C:Ck-circ-nonzero},
$v_p(C_k) \ge \frac{S^p(k)}{p-1}$,
with equality if and only if
$\overline{\bC_k^\circ}(a) \ne 0$.

Recall the polygon $\GL(d_1,d_2,p)$
from \cref{D:ANP}.
Consider the product
\begin{equation}
    \label{E:bCp}
    \bC^p(A) = \prod_{k \in V_p} \bC_k^\circ(A).
\end{equation}
Since
$\bC_k^\circ$ has nonzero reduction modulo $p$ for every $k\in V_p$,
$\bC^p(A)$ is a nonzero homogeneous polynomial as well,
and with coefficients in $\Q \cap \Z_p$.
Moreover, $\overline{\bC^p}(a) \ne 0$ if and only if 
$v_p(C_k) = \frac{S^p(k)}{p-1}$ for all $k \in V_p$.
In \cref{S:extension_theorem},
we will show that these are precisely the conditions needed to
achieve the generic Newton polygon for $L^*(f; s)$.

In fact, we can already prove more about $\bC^p$.
At the cost of asking $p$ to be slightly larger if necessary,
we can actually assert that the polynomial $\bC^p$
decomposes into a simple monomial component
that depends on $p$,
and an elaborate homogeneous component that doesn't,
in the sense illustrated by the following lemma:
\begin{lemma}
    \label{L:modular-invariant}
    Let 
    $r \in (\Z / D\Z)^\times$.
    Then there exists a nonzero 
    homogeneous polynomial $H^{(r)}(A) \in \Q[A]$
    such that, for any prime $p > \max\{2D, d_1^2, d_2^2\}$ 
    with $p \equiv r \pmod{D}$,
    $H^{(r)}\in(\Z_p\cap\Q)[A]$ and
    \[
        \bC^p(A) \equiv U_p(A) \cdot H^{(r)}(A) \pmod{p},
    \]
    where 
    \[
        U_p(A) = u_p A_{-d_1}^{e_{p,-}} A_{d_2}^{e_{p,+}},
        \qquad u_p \in \Z_p^\times \cap \Q,
        \quad e_{p,\pm} \in \Z_{\ge 0}.
    \]
\end{lemma}
\begin{proof}
    We first record a uniform form of the minimal-representation expansion.
    Let $(z)_v = z(z-1)\cdots(z-v+1)$ denote the falling factorial of $z$
    over $v$ iterations, with $(z)_0 = 1$.

    Let $\Delta$, $d$ be integers satisfying 
    $0 \le \Delta < d$.
    For nonnegative integers $v$ and $n$,
    let $\cT_{v}(n)$ denote the set of nonnegative integer tuples
    $\mathbf{u} = (u_1, \dots, u_v) \in \Z_{\ge 0}^v$ for which
    $\sum_{i=1}^v u_i i = n$,
    with the convention $\cT_0(0) = \{()\}$ and $\cT_0(n) = \varnothing$
    for positive $n$.
    We define
    \[
        P_{d, \Delta}(z; B_1, \dots, B_d)
        \coloneqq 
        \sum_{\mathbf{u} \in \cT_{d-1}(\Delta)}
        \frac{(z)_{|\mathbf{u}|}}{\prod_{i=1}^{d-1}u_{d-i}!}
        B_d^{d-1-|\mathbf{u}|}\prod_{i=1}^{d-1} B_i^{u_{d-i}}.
    \]
    Because $0 \le \Delta < d$, $P_{d, \Delta}$ is a homogeneous degree $d-1$
    polynomial in the variables $B_1, \dots, B_d$.

    Choose $N$ to be a nonzero integer for which
    $d-1 \le M(N) < p$,
    where we set $d = d_1$ if $N<0$
    and $d = d_2$ if $N>0$. 
    Then for any $\bt = (t_i) \in \Sol(N)$,
    the coefficient $\beta_{t_i}$ of $s^{t_i}$ in the Artin-Hasse exponential
    function $E_p(s)$ is equal to $(t_i)!^{-1}$.
    Thus, combining \cref{L:lower-term-bijection} 
    with
    \eqref{E:bGn-min}, we obtain
    \begin{equation}
        \label{E:bGN-factorization}
        \bG_N^\circ(A) = 
        \frac{A_{\varepsilon d}^{M(N)-d+1}}{M(N)!}
        P_{d, \Delta_N}(M(N); A_{\varepsilon 1}, \dots, A_{\varepsilon d}),
    \end{equation}
    where $\Delta_N = dM(N) - |N|$ and $\varepsilon = \sgn(N)$.

    Let $r \in (\Z / D\Z)^\times$,
    and for $i\in\{1,2\}$,
    let $r_i$ denote the image of $r$ in $(\Z / d_i\Z)^\times$.
    Let $p > \max\{2D, d_1^2, d_2^2\}$ 
    be a prime congruent to $r$ modulo $D$.
    By \cref{C:uniqueness-is-invariant},
    the set $V_p$ depends only on $r$,
    as do the unique minimizing intervals
    $\cI = \{-k_1, \dots, k_2\}$ for each $k \in V_p$.
    Thus we notate $V_p$ as $V_r$
    for the remainder of this proof,
    and fix such a $k$ and $\cI = \{-k_1, \dots, k_2\}$.
    For every nonzero $\ell \in \cI$,
    we set
    \[
        (\varepsilon_\ell, \delta_\ell, \kappa_\ell, \rho_\ell) = 
        \begin{cases*}
            (\phantom{-}1, d_2, k_2, r_2) & if $\ell > 0$ \\
            (-1, d_1, k_1, r_1) & if $\ell < 0$.
        \end{cases*}
    \]
    Now set
    \[
        e_\ell = [\kappa_\ell - \rho_\ell\cdot|\ell|]_{\delta_\ell}, \quad 
        m_\ell = \frac{e_\ell + p|\ell|-\kappa_\ell}{\delta_\ell}, \quad 
        c_\ell = \frac{e_\ell - \kappa_\ell}{\delta_\ell}.
    \]
    Note that 
    $e_\ell$ and $m_\ell$ are integers, $c_\ell \in \Z_p \cap \Q$,
    and $m_\ell \equiv c_\ell \pmod{p}$.
    Moreover, the bound on $p$, along with the inequality 
    $|\ell| \le \kappa_\ell \le \delta_\ell$, 
    gives $\delta_\ell-1 \le m_\ell < p$.
    
    Note that, as argued in \cref{R:SpI-mod-invariant},
    the set $\bSig^p(\cI)$ only depends on $r$. 
    We emphasize this by writing $\bSig^r(\cI)$ in place of 
    $\bSig^p(\cI)$.
    If $\sigma \in \bSig^r(\cI)$, then
    by \cref{P:small-one},
    $M(p\ell - \sigma\ell) = 
    \left\lceil \frac{p|\ell| - \kappa_\ell}{\delta_\ell}
        \right\rceil = m_\ell$.
    Now set 
    \[
        \Delta_{\ell, \sigma} = \delta_\ell m_\ell - |p\ell - \sigma\ell|
        = e_\ell - \kappa_\ell + \varepsilon_\ell 
        \sigma\ell \in \{0, \dots, \delta_\ell-1\}.
    \]
    Notice that $\Delta_{\ell, \sigma}$ does not depend on $p$,
    even though $\Delta_{p\ell - \sigma\ell} = \Delta_{\ell, \sigma}$.

    Since $\delta_\ell-1 \le m_\ell < p$,
    we can apply \eqref{E:bGN-factorization} to 
    \eqref{E:bCk-min} 
    and reduce modulo $p$
    to obtain
    \[
        \bC_k^\circ(A) \equiv U_{p, k}(A) H_k^{(r)}(A) \pmod{p},
    \]
    where
    \[
        U_{p, k}(A) = \prod_{\ell \in \cI\setminus \{0\}}
        \frac{A_{\varepsilon_\ell \delta_\ell}^{m_\ell - \delta_\ell + 1}}
            {m_\ell!}
    \]
    and 
    \[
        H_k^{(r)}(A) = 
        (-1)^k \sum_{\sigma \in \bSig^{r}(\cI)}
        \sgn(\sigma) \prod_{\ell \in \cI \setminus\{0\}}
        P_{\delta_\ell, \Delta_{\ell, \sigma}}(c_\ell, A_{\varepsilon_\ell 1}, 
        \dots, A_{\varepsilon_\ell \delta_\ell}).
    \]
    We see that $H_k^{(r)}(A)$ does not depend on $p$,
    and it is homogeneous of degree $k_1(d_1-1) + k_2(d_2-1)$,
    and has rational coefficients that are $p$-integral.
    The factor $U_{p, k}(A)$ is a monomial exclusively in $A_{-d_1}$
    and $A_{d_2}$ and the coefficient is a $p$-adic unit.
    Setting $U_p(A) = \prod_{k \in V_r} U_{p, k}(A)$
    and $H^{(r)}(A) = \prod_{k \in V_r} H_k^{(r)}(A)$
    implies the factorization congruence
    \[
        \bC^p(A) \equiv U_p(A) \cdot H^{(r)}(A) \pmod{p}.
    \]
    Note that this congruence itself implies that $H^{(r)}$
    is nonzero, since
    $\bC^p$ is nonzero modulo $p$ by \cref{C:Ck-circ-nonzero}.
    Therefore the lemma is proven.
\end{proof}

\begin{remark}[Cross-sign minimizing permutations]
    \label{R:cross-sign-minimizers}
    The positive and negative components of a Laurent polynomial
    do not, in general, contribute independently to the
    minimal terms of the Fredholm determinant.
    For example, let $d_1=d_2=3$ and $p=11$.
    For $k=3$, the unique minimizing index set is
    $\cI=\{-1,0,1\}$, with $S^{11}(3)=8$.
    Both the identity permutation and the transposition
    $\sigma=(-1\;\;1)$, fixing $0$, attain this minimum:
    \[
        M(-10)+M(10)=M(-12)+M(12)=8.
    \]
    In particular, the latter minimizing permutation
    does not preserve signs.

    The two permutations contribute with opposite signs
    to the minimal homogeneous component, giving
    \[
        \bC_3^\circ
        =-\bG_{-10}^\circ\bG_{10}^\circ
        +\bG_{-12}^\circ\bG_{12}^\circ
        =-P_-P_+
        +\frac{A_{-3}^4A_3^4}{(4!)^2},
    \]
    where
    \[
        P_\pm =
        \frac{A_{\pm3}^3A_{\pm1}}{3!}
        +\frac{A_{\pm3}^2A_{\pm2}^2}{(2!)^2}.
    \]
    Thus $\bC_3^\circ$ cannot be expressed as the product
    of a polynomial in the negative-indexed variables
    and a polynomial in the positive-indexed variables.

    In this example, $V_{11}=\{1,3,5\}$, with
    \[
        \bC_1^\circ=-1,
        \qquad
        \bC_5^\circ
        =-\frac{A_{-3}^{10}A_3^{10}}{(3!7!)^2}.
    \]
    Hence the genericity polynomial $\bC^{11}$
    likewise does not separate into independent
    positive and negative factors.
    This demonstrates why the two-pole case 
    we are studying in this paper does not
    follow simply by combining the corresponding
    one-pole constructions.
\end{remark}

\begin{example}
    We shall compute $H^{(r)} = H^{(r)}_1 H^{(r)}_3 H^{(r)}_5$ 
    for the example given in 
    \cref{R:cross-sign-minimizers}.
    For $k=1$,
    we have $\cI = \{0\}$, 
    so $H_1^{(r)} = -1$.
    
    For $k=3$, $\cI = \{-1, 0, 1\}$,
    we have $c_{-1} = c_1 = \frac{1}{3}$.
    There are two minimizing permutations
    $\sigma_1 = \id$ and $\sigma_2 = (-1\;\;1)$.
    They give $\Delta_{-1, \sigma_1} = \Delta_{1, \sigma_1} = 2$
    and $\Delta_{-1, \sigma_2} = \Delta_{1, \sigma_2} = 0$.
    Thus,
    \begin{align*}
        H_3^{(r)} &= -P_{3,2}(\frac{1}{3},A_{-1},A_{-2},A_{-3})
        P_{3,2}(\frac{1}{3},A_1,A_2,A_3) \\
                  &\quad + 
        P_{3,0}(\frac{1}{3}, A_{-1},A_{-2},A_{-3})
        P_{3,0}(\frac{1}{3}, A_1,A_2,A_3) \\
                  &= -\left(-\frac{1}{9}A_{-2}^2 + \frac{1}{3}A_{-1}A_{-3}\right)
                  \left(-\frac{1}{9}A_2^2 + \frac{1}{3}A_{1}A_{3}\right) \\
                  &\quad + A_{-3}^2 A_3^2.
    \end{align*}
    
    For $k=5$, the permutation $\sigma = (1\;\;2)(-1\;\;-2)$
    is the only minimizing permutation.
    We have $\Delta_{\ell, \sigma} = 0$ 
    for all $\ell \in \cI \setminus\{0\}$,
    and for any $c$, $P_{3, 0}(c, B_1, B_2, B_3) = B_3^2$.
    Thus
    \[
        H_5^{(r)} = -A_{-3}^4 A_3^4.
    \]
    This yields
    \[
        H^{(r)} = \frac{1}{81} A_{-3}^4A_3^4
        \left(81 A_{-3}^2 A_3^2 - 
        \left(3A_{-1}A_{-3}-A_{-2}^2\right)(3A_1A_3 - A_2^2)\right).
    \]
    
\end{example}

 \section{Genericity of the Fredholm determinant for field extensions}
\label{S:extension_theorem}
In this brief section,
we prove 
a proposition that effectively reduces the
generic condition to the base field case.

Recall the value $S^p(k)$ from \cref{SS:cycle-sums}.
Also recall the combinatorial cost polygon $\GL(d_1,d_2,p)$ 
from \cref{D:ANP}.
When $p > 2D$,
we have by \cref{C:kplus1} that
$k \in V_p$ if and only if $S^p(\cI) = S^p(k)$
for a unique $\cI \in \binom{\Z}{k}$.

We argued in \cref{SS:pind-factor}
that 
for sufficiently large $p$,
$\GL(d_1,d_2,p)$ is the generic polygon
for the Fredholm determinant $\det(1 - \cG s)$
defined in \cref{SS:dwork-operator}.
We wish to establish $\GL(d_1,d_2,p)$ to be the generic Newton polygon
for the $L$-function $L^*(f; s)$ as well.
But \eqref{E:L*-G-tau-b} connects the Newton polygon of $L^*(f; s)$
with the product matrix 
$[\cG, \tau, b] = \cG^{\tau^{b-1}} \cdots \cG^\tau \cG$
given in \eqref{E:matGb}, where $\tau$ is the lift of the Frobenius.
Therefore we must connect the Newton polygon of the Fredholm determinant 
of $[\cG, \tau, b]$ with that of $\cG$,
and then with $\GL(d_1,d_2,p)$.
We establish this connection in the proposition below.
\begin{proposition}
    \label{P:generic-np-extension}
    Assume $p > 2D$.
    Let $q = p^b$, and let $\tau$
    denote the lift of the Frobenius automorphism
    on $\F_q$ to
    $\Q_q(\zeta_p)$, extended to fix $\zeta_p$.
    Let $a = (a_i)_{i \in I_0} \in \F_q^{I_0}$,
    where $I_0 = \{-d_1, \dots, d_2\} \setminus \{0\}$,
    and suppose $a_{-d_1}a_{d_2} \ne 0$.
    For every $N \in \Z$, set 
    $G_N = {\bG}_N(\wh{a})$,
    and define $\cG = (G_{pi-j})_{i,j \in \Z}$
    and $[\cG, \tau, b] = \cG^{\tau^{b-1}} \cdots \cG^\tau \cG$.

    Then 
    \[
        \NP_q(\det(1 - [\cG, \tau, b]s))\big|_{[0,d_1+d_2]} \ge 
        \GL(d_1,d_2,p).
    \]
    Moreover,
    equality holds if and only if 
    $v_p(C_k) = \frac{S^p(k)}{p-1}$ for all $k \in V_p$.
\end{proposition}
\begin{proof}
    A Cauchy-Binet--type argument produces the following result.
    \begin{lemma}[Cauchy-Binet]
        \label{L:cauchy-binet}
        For every $m \in \Z / b\Z$, let 
        $A_m$ be a nuclear matrix expressed with a formal basis
        $\{e_i\}_{i \in \Z}$,
        and let $P = A_{b-1}A_{b-2}\cdots A_0$,
        and let 
        \[
            \det(1 - Ps) = 1 + P_1s + P_2s^2 + \dots.
        \]
        Then for every $k \ge 1$,
        \[
            P_k = (-1)^k \sum_{\cI_0, \dots, \cI_{b-1} \in \binom{\Z}{k}}
            \prod_{m \in \Z / b\Z} \det((A_m)_{\cI_m, \cI_{m-1}}),
        \]
        where $(A_m)_{\cI_m, \cI_{m-1}}$ denotes the submatrix of $A_m$
        with rows indexed by $\cI_m$ and columns indexed by $\cI_{m-1}$.
    \end{lemma}
    The proof involves a routine argument using finite truncations,
    leveraging the convergence of the Fredholm expansion of the product $P$.
    We omit the details here.

    Assume $1 \le k \le d_1+d_2$, and 
    let $D_k$ denote the coefficient of $s^k$ in 
    $\det(1-[\cG, \tau, b]s)$.
    By \cref{L:cauchy-binet}, we have that
    \[
        v_p(D_k) \ge \min_{\cI_0, \dots, \cI_{b-1}} 
        \sum_{m \in \Z/ b\Z} v_p(\det(\cG^{\tau^m}_{\cI_m, \cI_{m-1}})),
    \]
    with equality if there is a unique tuple $(\cI_m)_{m \in \Z/b\Z}$
    that achieves the minimum.
    Note that in general, 
    $v_p(\det(N)) = v_p(\det(N^{\tau^{i}}))$ for any finite square matrix
    $N$, since $v_p$ is preserved under automorphism, 
    so we restrict our attention to $k \times k$ minors of $\cG$
    in the argument that follows.
    For $\cI, \mathcal{J} \in \binom{\Z}{k}$,
    let $\Delta[\cI, \mathcal{J}] = v_p(\det(\cG_{\cI, \mathcal{J}}))$.

    Choose any tuple $(\cI_m) \in \binom{\Z}{k}^{\Z / b\Z}$.
    Then 
    \[
        \det(\cG_{\cI_m, \cI_{m-1}}) = 
        \sum_{\sigma}\sgn(\sigma) \prod_{\ell \in \cI_m}
        G_{p\ell - \sigma\ell},
    \]
    where $\sum_{\sigma}$ ranges over all bijections 
    $\sigma : \cI_m \to \cI_{m-1}$.
    Thus, we have by an application of \cref{P:bG-expression} that 
    \begin{equation}
        \label{E:fred-ineq}
        \Delta[\cI_m, \cI_{m-1}]
        \ge \min_\sigma \sum_{\ell \in \cI_{m}} 
        \frac{M(p\ell-\sigma\ell)}{p-1}
        = \frac{S^p(\cI_{m}, \cI_{m-1})}{p-1},
    \end{equation}
    where $S^p(\cI_m, \cI_{m-1})$ is equal to
    the minimum value of $\sum_{\ell \in \cI_m} M(p\ell - \sigma\ell)$
    as $\sigma$ varies over all bijections $\cI_m \to \cI_{m-1}$.
    After combining the above inequalities, 
    we have by \cref{P:minimal-bijections} that
    $v_p(D_k) \ge (bS^p(k))/(p-1)$, or, after normalization,
    \[
        v_q(D_k) 
        \ge \frac{S^p(k)}{p-1}.
    \]
    These coefficient bounds suffice to establish the claimed inequality 
    for the full Newton polygon --- that is,
    $\NP_q(\det(1 - [\cG, \tau, b]s))|_{[0,d_1+d_2]} \ge \GL(d_1,d_2,p)$.
    Indeed, by \cref{prop:L-fredholm-NP} and \eqref{E:L*-G-tau-b},
    the Newton polygon of $\det(1-[\cG,\tau,b]s)$ agrees on $[0,d_1+d_2]$
    with that of $L^*(f;s)$,
    which has slopes in $[0,1]$ and endpoint $(d_1+d_2,(d_1+d_2)/2)$.
    As we verify below, $\GL(d_1,d_2,p)$ has the same endpoint, 
    and its slopes are at most $1$ by \cref{C:kplus1}. 
    Thus, coefficients of degree greater than $d_1+d_2$
    cannot lower the Newton polygon beneath $\GL(d_1,d_2,p)$
    on the specified interval,
    and the preceding inequalities establish the claim.

    Let $k \in V_p$.
    By \cref{C:kplus1}, there must be a unique $\cI \in \binom{\Z}{k}$
    for which $S^p(\cI) = S^p(k)$.
    Then 
    by \eqref{E:fred-ineq} and \cref{P:minimal-bijections},
    \[
        \sum_{m \in \Z / b\Z} \Delta[\cI_m, \cI_{m-1}]
        > \frac{bS^p(k)}{p-1}
    \]
    if $\cI_m \ne \cI$ for some $m \in \Z / b\Z$,
    and in particular,
    \[
        \Delta[\cI', \cI'] > \frac{S^p(k)}{p-1}
    \]
    if $\cI' \ne \cI$.
    Since 
    \[
        v_p(D_k) \ge \min_{\cI_0, \dots, \cI_{b-1}} 
        \sum_m \Delta[\cI_m, \cI_{m-1}] 
        \quad \mbox{and} \quad
        v_p(C_k) \ge \min_{\cI'} \Delta[\cI', \cI']
    \]
    with equalities if the minima are achieved uniquely,
    we have
    \[
        v_q(D_k) = \frac{S^p(k)}{p-1} \iff \Delta[\cI, \cI] = \frac{S^p(k)}{p-1}
        \iff v_p(C_k) = \frac{S^p(k)}{p-1}.
    \]
    If
    \[
        \NP_q(\det(1 - [\cG, \tau, b]s))|_{[0, d_1+d_2]} = \GL(d_1,d_2,p),
    \]
    then $v_q(D_k) = \frac{S^p(k)}{p-1}$,
    thus $v_p(C_k) = \frac{S^p(k)}{p-1}$,
    for all $k \in V_p$.
    Conversely, if $v_p(C_k) = \frac{S^p(k)}{p-1}$ 
    for all interior vertices $k \in V_p$,
    we also have that $\NP_q(\det(1 - [\cG, \tau, b]s))$
    and $\GL(d_1,d_2,p)$ agree on all interior vertices $k \in V_p$.
    It remains to show that these two polygons agree on the endpoints as well.
    For $k=0$, this corresponds to the origin $(0,0)$ on both polygons.
    For $k=d_1+d_2$,
    we have $S^p(d_1 + d_2) = \frac{(p-1)(d_1+d_2)}{2}$.
    Indeed, choose the minimizing set $\cI = \{-d_1, \dots, d_2-1\}$.
    Choose a sign-preserving permutation $\sigma$ on $\cI$ that satisfies 
    $p\ell \equiv \sigma\ell \pmod{d_1}$ for all $\ell < 0$
    and $p\ell \equiv \sigma\ell \pmod{d_2}$ for all $\ell > 0$,
    and $\sigma 0 = 0$.
    Such a permutation exists because $p$ is coprime to $d_1d_2$.
    Choosing such a $\sigma$ makes all relevant quotients integral, and
    we have
    \[
        S^p(k) \le S^p(\cI) \le \frac{p-1}{d_1}(1 + \dots + d_1)
        + \frac{p-1}{d_2}(1 + \dots + d_2-1) = \frac{(p-1)(d_1+d_2)}{2}.
    \]
    This is also a lower bound for $S^p(k)$,
    since $W_{d_1+d_2} = w(\cI) = \frac{d_1 + d_2}{2}$,
    and $S^p(k) \ge (p-1)W_{d_1+d_2}$ by \cref{L:weight-properties} (5).
    Thus $\GL(d_1,d_2,p)(d_1+d_2) = \frac{d_1+d_2}{2}$.

    On the other hand, by geometry, the $L$-function
    $L^*(f; s)$ always has an endpoint of 
    $\left(d_1+d_2, \frac{d_1+d_2}{2}\right)$ on its $q$-adic Newton polygon.
    Indeed, it is a degree $d_1+d_2$ polynomial with slopes coming in 
    pairs of $\lambda$ and $1-\lambda$.
    By \eqref{E:L*-G-tau-b}, the same endpoint lies on 
    $\NP_q(\det(1 - [\cG, \tau, b]s))|_{[0,d_1+d_2]}$.
\end{proof}
 \section{Main results}
\label{S:independence}

Here, we combine the results of \cref{S:ck_formula,S:extension_theorem}
to establish the generic condition for the $L$-function.
Our main result will then follow.

\cref{S:extension_theorem} established the condition for genericity
for the Newton polygon of the 
Fredholm determinant $\det(1 - [\cG, \tau, b]s)$ on the interval 
$[0, d_1+d_2]$.
The matrix $[\cG, \tau, b]$ is the matrix form of the Dwork operator
associated with a Laurent polynomial $f$,
and so by \eqref{E:L*-G-tau-b}, \cref{P:generic-np-extension}
also establishes the generic condition for $L^*(f; s)$.
Namely, this proposition reduces the genericity of $L^*(f; s)$
to showing that the coefficients $C_k$
of the Fredholm determinant $\det(1 - \cG s) = 1 + C_1s + \dots$
align with the combinatorial cost polygon $\GL(d_1,d_2,p)$
on its internal vertices.

\subsection{Proof of \Cref{T:ordinary-curves}}

Recall the expression $\bC^p$ given by
\eqref{E:bCp}.
This was proven to be a nonzero homogeneous polynomial
in $(\Z_p \cap \Q)[A]$.
We are ready now to prove that $\bC^p$
describes the generic condition for $L^*(f; s)$.
\begin{proposition}
    \label{T:generic-alignment}
    Let $
        f = \sum_{i \in I_0} a_i x^i \in \F_{q}[x^{\pm1}]$,
 where $a_{-d_1} a_{d_2} \ne 0.$
    Suppose $p > 2D$,
    and let $\overline{\bC^p}$ denote the
    reduction of $\bC^p$ in $\overline{\F}_p$.
    Then
    \[
        \NP_q(L^*(f; s)) \ge \GL(d_1,d_2,p),
    \]
    with equality if and only if
    $\overline{\bC^p}(a) \ne 0$.
\end{proposition}
\begin{proof}
    By combining
    \eqref{E:L*-G-tau-b},
    and \cref{P:generic-np-extension},
    we immediately obtain
    \[
        \NP_q(L^*(f; s)) \ge \GL(d_1,d_2,p),
    \]
    with equality if and only if 
    $v_p(C_k) = \frac{S^p(k)}{p-1}$ for every $k \in V_p$.
    Thus it remains to show that this is equivalent to
    $\overline{\bC^p}(a) \ne 0$.
    Choose $k \in V_p$.
    The coefficient $C_k$ of $s^k$ in $\det(1 - \cG s)$
    is given by the evaluation of $\bC_k(\wh{a})$,
    where $\wh{a}$ denotes the Teichm\"uller lift of 
    the $a_i$.
    It now follows by \cref{C:Ck-circ-nonzero} that
    \[
        v_p(C_k) = \frac{S^p(k)}{p-1} 
        \iff
        \bC_k^\circ(\wh{a}) \not\equiv 0 \pmod{\km_{\overline{\Q}_p}}
        \iff 
        \overline{\bC_k^\circ}(a) \ne 0.
    \]
    Thus $v_p(C_k) = \frac{S^p(k)}{p-1}$ for every $k \in V_p$
    if and only if $\overline{\bC^p}(a) = \prod_{k \in V_p}
    \overline{\bC_k^\circ}(a) \ne 0$.
\end{proof}

We now proceed to proving our main results.
\begin{theorem}
    \label{T:open-subset}
    For any $a \in \A^I_\Q(\overline{\Q})$,
    let $f_a$ denote the Laurent polynomial $\sum_{i \in I} a_i x^i$
    associated with $a$.
    Then there is a homogeneous Zariski-open and dense subset 
    $\mathcal{U} \subseteq \A^I_{\Q}$
    defined over $\Q$ for which the following are equivalent
    for all $a \in \A^I_\Q(\overline{\Q})$
    with $a_{-d_1}a_{d_2} \ne 0$:
    \begin{enumerate}[(1)]
        \item $a \in \mathcal{U}(\overline{\Q})$,
        \item There is a constant $N_a$ such that,
            for every prime number $p > N_a$,
            every number field $K / \Q(a_i \mid i \in I)$,
            and every prime ideal $\kp$ in $K$ above $p$,
            the reduction $\overline{f}_a = f_a \pmod{\kp}$
            is defined, and for every $\ell$ that divides 
            $[k_\kp : \F_p]$ and every nontrivial 
            additive character 
            $\chi_\ell : \F_{p^\ell} \to \overline{\Q}_p^\times$,
            the $L$-function $L^*(\overline{f}_a, \chi_\ell; s)$
            has generic Newton polygon.
    \end{enumerate}
\end{theorem}
\begin{proof}
    Let
    \[
        H(A) = A_{-d_1}A_{d_2} \prod_{r \in (\Z / D\Z)^\times} H^{(r)}(A) 
        \quad \mbox{ and } \quad \mathcal{U} = D(H).
    \]
    As defined, this is a polynomial in the variables indexed by $I_0$.    Extend it to a polynomial with variables indexed by $I$
    in the obvious way,
    so that the constant $A_0$ is free.
    Since $H$ is nonzero and homogeneous, and defined over $\Q$,
    $\mathcal{U}$ is a dense homogeneous open set defined over $\Q$.
    
    Let $a \in \mathcal{U}(\overline{\Q})$.
    Then $H^{(r)}(a) \ne 0$ for every $r$.
    Let $E = \Q(a_i \mid i \in I)$,
    and note that the following sets are finite:
    \begin{enumerate}[i.]
        \item The set of rational primes below places of $E$ for which 
            $a_i$ is non-integral for some $i \in I$,
        \item The set of rational primes below places of $E$ for which 
            $a_{-d_1}$ or $a_{d_2}$ is not a unit,
        \item The set of rational primes below places of $E$ for which
            $H^{(r)}(a)$ is not a unit for some $r \in (\Z / D\Z)^\times$,
        \item The set of rational primes $p$ for which 
            $p \le \max \{ 2D, d_1^2, d_2^2\}$.
    \end{enumerate}
    Let $N_a$ denote the largest prime contained in the union of these sets,
    and let $p > N_a$.
    Let $K / E$ be a number field, and let $\kp$ be a prime ideal above $p$.
    By first applying the shift 
    $\overline{f}_a \mapsto \overline{f}_a - \overline{a}_0$
    as described in \cref{r:assume-a0-zero}, we may assume that 
    $\overline{a}_0 = 0$.
    Then 
    $\overline{H^{(r)}}(\overline{a})$ in $k_\kp$
    is defined and nonzero for all $r \in (\Z / D\Z)^\times$.
    Take $r$ to be the the residue class of $p\bmod D$, 
    \cref{L:modular-invariant} gives that
    \[
        \overline{\bC^p}(\overline{a}) =
        \overline{U_p}(\overline{a}) \overline{H^{(r)}}(\overline{a}) \ne 0.
    \]
    We thus have by \cref{T:generic-alignment} that
    $L^*(\overline{f}_a; s)$ has generic Newton polygon.
    
    Now, if $\chi_\ell : \F_{p^\ell} \to \overline{\Q}_p^\times$
    is a nontrivial additive character,
    then by \cref{P:rank-l-reduction}, there is an 
    $\alpha \in \F_{p^\ell}^\times$
    for which $L^*(\overline{f}_a, \chi_\ell; s) 
    = L^*(\alpha \overline{f}_a, s)$.
    Since $\overline{\bC^p}(\bar{a}) \ne 0$
    and $\bC^p$ is homogeneous, it follows that 
    $\overline{\bC^p}(\alpha \overline{a}) \ne 0$ in $k_{\kp}$,
    and by \cref{T:generic-alignment} we conclude that 
    $L^*(\alpha \overline{f}_a, s)$ has generic Newton polygon.

    Conversely, suppose $a \not\in \mathcal{U}(\overline{\Q})$.
    Then since $a_{-d_1}a_{d_2} \ne 0$,
    we will have $a_{-d_1}a_{d_2} \not\equiv 0 \pmod{p}$
    for sufficiently large $p$.
    Therefore $H^{(r)}(a) = 0$ for some $r \in (\Z / D\Z)^\times$.
    By Dirichlet's theorem, there are infinitely many primes $p$
    in this residue class.
    The sets i., ii., and iv. above are still finite sets;
    thus there are also infinitely many primes $p'$
    outside of these sets with the same congruence class.
    For each such $p'$, choose a prime ideal $\kp'$ above $p'$.
    Then since $H^{(r)}(a) = 0$, we have by \cref{L:modular-invariant}
    that
    $\overline{\bC^{p'}}(\overline{a}) = 0$.
    By \cref{T:generic-alignment}, the $L$-function
    $L^*(\overline{f}_{a}; s)$ will not be generic.
    Therefore there is no such constant $N_a$.
\end{proof}
\begin{theorem}[\cref{T:ordinary-curves}]
    With $I$, $f_a$, and $\mathcal{U}$ as above,
    suppose $a_{-d_1}a_{d_2} \ne 0$.
    Then the following are equivalent: 
    \begin{enumerate}[(1)]
        \item $a \in \mathcal{U}(\overline{\Q})$
        \item There is a constant $N_a$ such that, for every prime number
            $p > N_a$, every number field $K / \Q(a_i \mid i \in I)$,
            every prime ideal $\kp$ in $K$ above $p$, and every $\ell$
            that divides $[k_\kp : \F_p]$,
            the reduction $\overline{f}_a = f_a \pmod{\kp}$
            is defined, and the curve $X_{\overline{f}_a, \ell}$
            is generic.
    \end{enumerate}
\end{theorem}
\begin{proof}
    Let $a \in \mathcal{U}(\overline{\Q})$.
   Let $N_a$ be a constant for which condition (2) of \cref{T:open-subset} holds.
    Choose $p > N_a$, $K / \Q(a_i \mid i \in I)$,
    $\kp$ above $p$, and $\ell$ dividing $[k_\kp : \F_p]$.
    Note that \cref{T:open-subset}
    implies that $L^*(\overline{f}_a; s)$ has generic Newton polygon.
    By the homogeneity of the reduction $\bC^p$ mod $\kp$,
    $L^*(\lambda \overline{f}_a; s)$ 
    is generic for every $\lambda \in \F_{p^\ell}^\times$.
    By \cref{P:zeta-L-connection},
    the numerator of the zeta function for 
    $X_{\overline{f}_a, \ell}$ will be
    \[
        \prod_{\lambda \in \F_{p^\ell}^\times} 
        L^*(\lambda \overline{f}_a; s).
    \]
    The Newton slopes of this numerator are equal to the multiset union 
    of the slopes of each factor $L^*(\lambda \overline{f}_a; s)$.
    Since every factor is bounded below by the same generic Newton polygon,
    and all generic conditions are simultaneously obtained,
    it follows that $X_{\overline{f}_a, \ell}$
    is generic if and only if $L^*(\lambda \overline{f}_a; s)$
    has generic Newton polygon for every $\lambda \in \F_{p^\ell}^\times$.
    This proves the forward direction.

    Conversely, 
If $a\notin\mathcal U(\overline{\mathbb Q})$, 
the converse argument in \Cref{T:open-subset}
produces arbitrarily large primes and reductions for which 
$L^*(\overline f_a;s)$ is nongeneric. 
For $\ell=1$, this is one of the factors in the zeta numerator of $X_{\overline f_a,1}$. Since every factor is bounded below by the generic polygon, the curve is nongeneric. Thus no constant $N_a$ can satisfy condition (2).
\end{proof}

\subsection{Proof of \Cref{cor:introduction-T-adic}}
%C:kplus1

Recall notations from \cref{SS:T-adic-intro}.
In \Cref{P:lln-result} below, we establish an interpolation criterion under the hypothesis $p\ge2D$, using the method of \cite{ref:liu2016}; 
its global consequence is \Cref{T:T-adic}.
For $\Delta=[-d_1,d_2]$, let $p_\Delta$ denote the infinite arithmetic
polygon of \cite[Definition 1.6]{ref:liu2016}.

\begin{proposition}[Interpolation criterion]
    \label{P:lln-result}
    Suppose $p\ge 2D$, and let $q=p^b$, $h=d_1+d_2$.
    Then there is a nonzero polynomial $H(y) \in \F_p[y_u \mid u \in I]$
    such that, for every 
     $f = \sum_{i \in I} a_i x^i \in \F_q[x^{\pm 1}]$
    with $a_{-d_1}a_{d_2} \ne 0$,
    the following are equivalent:
    \begin{enumerate}[(1)]
        \item $H(a) \ne 0$, where $a = (a_i)_{i \in I}$,
        \item $L^*(f; T, s)|_{T=\pi_1}$ has generic Newton polygon,
        \item $C^*(f; T, s)|_{T=\pi_1}$ has generic Newton polygon,
        \item $C^*(f; T, s)|_{T=\pi_m}$ has generic Newton polygon
            for every $m\ge 1$,
        \item $L^*(f; T, s)|_{T=\pi_m}$ has generic Newton polygon
            for every $m\ge 1$.
    \end{enumerate}
    Here the specialized polygons are taken with respect to the
    $\pi_m$-adic valuation. Their generic polygons are $b\,p_\Delta$
    for the $C$-functions and
    $b\,p_\Delta|_{[0,p^{m-1}h]}$ for the $L$-functions.

    Each of these equivalent conditions also implies
    \[
        \NP_T(C^*(f;T,s))=b\,p_\Delta,
    \]
    so that the $T$-adic $C$-function has generic Newton polygon.
\end{proposition}

\begin{proof}
    First suppose $p>2D$. We sharpen the arithmetic estimate underlying
    the interpolation argument of \cite{ref:liu2016}.
    For a positive integer $d$ dividing $D$, put
    \[
        r_d(n)=1+\frac{n}{d}-\left\lceil\frac{n}{d}\right\rceil,
        \qquad
        s_d(n)=r_d(pn),
        \qquad
        g_d(n)=\left\lceil\frac{pn}{d}\right\rceil
               -\left\lceil\frac{n}{d}\right\rceil.
    \]
    Thus $r_d(n),s_d(n)\in(0,1]$, and
    \[
        \left\lceil\frac{pi-j}{d}\right\rceil
        =
        \left\lceil\frac{pi}{d}\right\rceil
        -\left\lceil\frac{j}{d}\right\rceil
        +\mathbf{1}_{\{r_d(j)<s_d(i)\}}.
    \]
    If $A\subseteq\Z_{\ge0}\times\{1,\ldots,b\}$ is finite,
    $n(a)$ denotes the first coordinate of $a$, and $\tau$ is a
    permutation of $A$, define
    \[
        \Phi_d(A,\tau)=
        \sum_{a\in A}
        \left(g_d(n(a))+
        \mathbf{1}_{\{r_d(n(\tau a))<s_d(n(a))\}}\right).
    \]

    Suppose $A$ contains a labeled copy of $n$ and omits a labeled
    copy of $m<n$. Replace this copy of $n$ by the missing copy of $m$,
    and transport $\tau$ along the replacement bijection to obtain
    $(A',\tau')$. We claim that
    \[
        \Phi_d(A,\tau)\ge \Phi_d(A',\tau')+1.
    \]
    Every adjacent increment of $g_d$ is at least $2$: if $d\nmid k$,
    the increment is at least $\lfloor p/d\rfloor\ge2$, whereas if
    $d\mid k$, it is $\lceil p/d\rceil-1\ge2$.
    The replacement changes at most two indicators, one for the
    outgoing arrow and one for the incoming arrow.
    If the replaced element is fixed, there is only one indicator,
    and the claim follows. If $g_d(n)-g_d(m)\ge3$, the possible
    increase of the indicator sum is at most $2$, which also proves
    the claim.
    In the remaining case $g_d(n)-g_d(m)=2$, we have $n=m+1$.
    If $d\mid m$, then $r_d(m)=1\ge r_d(n)$, so the incoming
    indicator cannot increase.
    If $d\nmid m$, then
    $\lceil pn/d\rceil-\lceil pm/d\rceil=2$, whence
    $s_d(n)-s_d(m)=p/d-2>0$; the outgoing indicator cannot increase.
    Thus in either case the indicator sum increases by at most $1$,
    while the $g_d$-sum decreases by $2$. This proves the claim.

    Repeated replacements reduce $A$ to an initial labeled set.
    On such a set, the residue-count argument of
    \cite[Lemma 3.3]{ref:liu2016} uses only $(p,d)=1$.
    Consequently their one-sided arithmetic estimate, including
    strictness for noninitial sets, holds for $p>2d$.
    The same statement for positive indices follows by adjoining
    all $b$ copies of zero and fixing them.
    For clarity, the arithmetic slopes on this side are
    \[
        \omega_d(n)=
        \left\lceil\frac{(p-1)n}{d}\right\rceil-\epsilon_d(n),
        \qquad n\ge1,
    \]
    where $\epsilon_d(n)=1$ if there is an integer
    $1\le u<[n]_d$ with $pu\equiv n\pmod d$, and
    $\epsilon_d(n)=0$ otherwise.
    The slopes of $p_\Delta$ are the multiset consisting of $0$
    and the $\omega_{d_1}(n),\omega_{d_2}(n)$ for $n\ge1$.

    To pass to two signs, let
    $A\subseteq\Z\times\{1,\ldots,b\}$ have cardinality $bk$.
    Repair $\tau$ to a permutation $\tau_0$ which fixes the zero
    elements, preserves each nonzero sign, and retains every
    original arrow whose endpoints have the same nonzero sign.
    Complete the remaining arrows within each sign.
    For a nonzero source with first coordinate $i$, the original
    cost is at least
    \[
        \left\lceil
        \frac{p|i|-|n(\tau_0 a)|}{d(i)}
        \right\rceil,
        \qquad
        d(i)=
        \begin{cases}
            d_1,&i<0,\\
            d_2,&i>0.
        \end{cases}
    \]
    Indeed, this follows from the corresponding signed ceiling
    when the original arrow preserves sign; otherwise the original
    cost is at least $\lceil p|i|/d(i)\rceil$.
    The zero-source costs are nonnegative.
    Applying the one-sided estimate to each sign, and then selecting
    the $bk$ smallest slopes from $b$ copies of the arithmetic slope
    multiset, gives
    \[
        \sum_{a\in A}M\bigl(pn(a)-n(\tau a)\bigr)
        \ge b\,p_\Delta(k).
    \]
    This holds for every $k\ge1$. The $T$-adic Dwork trace formula,
    the semilinear-to-linear determinant comparison, and the Fredholm
    expansion, as in \cite[Theorems 2.1--2.2 and the proof of
    Theorem 4.1]{ref:liu2016}, therefore give
    \begin{equation}
        \label{E:T-adic-bound-2D}
        \NP_T(C^*(f;T,s))\ge b\,p_\Delta.
    \end{equation}

    The slope calculation in the proof of \cref{C:kplus1} identifies
    the first $h=d_1+d_2$ slopes and gives
    \[
        p_\Delta|_{[0,h]}=(p-1)\GL(d_1,d_2,p).
    \]
    Moreover,
    $\omega_d(n+d)=\omega_d(n)+(p-1)$.
    Thus the full slope multiset of $p_\Delta$ is obtained by
    repeating its first $h$ slopes with shifts $j(p-1)$,
    $j\ge0$.

    A constant term does not affect any of the polygons in the
    statement. More precisely, for $f_0=f-a_0$ and
    $c=\Tr_{\Q_q/\Q_p}(\widehat a_0)$, the change from $f_0$ to $f$
    replaces $s$ by $(1+T)^c s$ in both $C^*$ and $L^*$.
    This multiplier, and each of its specializations, is a unit.
    We may therefore assume $a_0=0$.
    By \cref{T:generic-alignment}, we may take
    $H=\overline{\bC^p}$, extended to be independent of $a_0$.
    It is nonzero, and
    \[
        H(a)\ne0
        \iff
        \NP_{\pi_1}\bigl(L^*(f;T,s)|_{T=\pi_1}\bigr)
        =b\,p_\Delta|_{[0,h]}.
    \]
    The product identity
    \[
        C^*(f;\pi_1,s)=
        \prod_{j\ge0}L^*(f;\pi_1,q^js),
    \]
    where specialization is indicated by writing $\pi_1$ in place of
    $T$, and the slope periodicity show that this is equivalent to
    $\NP_{\pi_1}(C^*(f;\pi_1,s))=b\,p_\Delta$.

    Write $C^*(f;T,s)=\sum_{i\ge0}c_i(T)s^i$.
    Its coefficients belong to $\Z_p[[T]]$, so specialization gives
    \[
        \NP_{\pi_m}(C^*(f;\pi_m,s))
        \ge \NP_T(C^*(f;T,s))\ge b\,p_\Delta.
    \]
    Equality at $\pi_1$ forces equality throughout this chain.
    At every vertex $(i,e)$, it also forces
    $c_i(T)=T^e u_i(T)$ with $u_i(T)\in\Z_p[[T]]^\times$.
    Such a unit remains a unit after every specialization
    $T=\pi_m$. Hence equality holds for every specialized
    $C$-function, and for the $T$-adic $C$-function.
    Conversely, equality at any specialization forces the same
    unit condition, and hence equality at $\pi_1$.

    Finally,
    \[
        C^*(f;\pi_m,s)=
        \prod_{j\ge0}L^*(f;\pi_m,q^js),
        \qquad
        \ord_{\pi_m}(q)=bp^{m-1}(p-1).
    \]
    The specialized $L$-function has degree $p^{m-1}h$.
    Comparing slopes in this product, as in
    \cite[Lemma 6.3]{ref:liu2016}, shows that the $C$-function
    attains $b\,p_\Delta$ if and only if the $L$-function attains
    its restriction to $[0,p^{m-1}h]$.
    Since the slope multiset of the $C$-function contains that of
    its $L$-function, the specialized $C$-function lower bound also
    bounds the $L$-function polygon below by $b\,p_\Delta$ on this
    interval.
    These are their generic polygons, since the preceding lower
    bounds are attained when $H(a)\ne0$.
    This proves the five equivalences and the $T$-adic consequence
    for $p>2D$.

    If $p=2D$, then $p=2$ and $d_1=d_2=D=1$.
    In this case take $H=1$.
    The classical $L$-function is the degree-two zeta numerator of
    the associated genus-one Artin--Schreier curve. Its linear
    coefficient is a sum of $q-1$ numbers in $\{1,-1\}$ and is
    therefore odd, while its leading coefficient is $q$.
    Its $q$-adic slopes are consequently $0,1$.
    When $D=1$, the general $T$-adic Hodge bound of
    \cite{ref:LiuWan2009} is exactly $b\,p_\Delta$.
    The same product and unit-coefficient arguments now give all
    the assertions in this remaining case.
\end{proof}

\begin{theorem}[\Cref{cor:introduction-T-adic}]
    \label{T:T-adic}
    Let $I$, $f_a$, and $\mathcal{U}$ be as in 
    \cref{T:open-subset}, and suppose $a_{-d_1}a_{d_2} \ne 0$.
    Then the following are equivalent:
    \begin{enumerate}[(1)]
        \item $a \in \mathcal{U}(\overline{\Q})$,
        \item There is a constant $N_a$ such that,
            for every prime number $p > N_a$,
            every number field $K / \Q(a_i \mid i \in I)$,
            every prime ideal $\kp$ in $K$ above $p$,
            and every $m \ge 1$,
            $C^*(\overline{f}_a; T, s)$, 
            $C^*(\overline{f}_a; T, s)|_{T = \pi_m}$,
            and $L^*(\overline{f}_a; T, s)|_{T = \pi_m}$
            all have generic Newton polygon.
    \end{enumerate}
\end{theorem}
\begin{proof}
    Suppose $a\in\mathcal{U}(\overline{\Q})$.
    By \cref{T:open-subset}, after enlarging $N_a$ so that
    $N_a\ge2D$, the classical $L$-function
    $L^*(\overline f_a;s)$ is generic for every prime and reduction
    in the statement.
    Equation \eqref{E:L-Tadic} and \cref{P:lln-result} then give
    genericity of all the specialized $C$- and $L$-functions,
    together with genericity of the $T$-adic $C$-function.

    Conversely, condition (2), with $m=1$, includes genericity of
    $L^*(\overline f_a;s)$ for all sufficiently large primes, by
    \eqref{E:L-Tadic}.
    The converse argument in the proof of \cref{T:open-subset}
    therefore gives $a\in\mathcal{U}(\overline{\Q})$.
\end{proof}

\medskip
\bibliographystyle{amsplain}
\providecommand{\bysame}{\leavevmode\hbox to3em{\hrulefill}\thinspace}
\providecommand{\MR}{\relax\ifhmode\unskip\space\fi MR }
% \MRhref is called by the amsart/book/proc definition of \MR.
\providecommand{\MRhref}[2]{%
  \href{http://www.ams.org/mathscinet-getitem?mr=#1}{#2}
}
\providecommand{\href}[2]{#2}

\end{document}